  \documentclass[12pt,final]{siamltex}

  \usepackage{times}
  \usepackage{amssymb}

\newcommand{\Rset}{\mathbb{R}}
\newcommand{\Zset}{\mathbb{Z}}
\newcommand{\Qset}{\mathbb{Q}}
\newcommand{\mfa}{\mathfrak{a}}
\newcommand{\mfb}{\mathfrak{b}}

\newcommand{\mft}{\mathfrak{t}}
\newcommand{\mfw}{\mathfrak{w}}

\newcommand{\bmatrix}[1]{\left[\matrix{#1}\right]}
\newtheorem{example}[theorem]{Example}
\newtheorem{note}[theorem]{Note}

\newcommand{\QEDBBOX}{\hspace*{\fill}\rule{0.4em}{0.4em}\par\bigskip}
\newcommand{\adj}{\mathop{\mathrm{adj}}\nolimits}

\renewcommand{\Bigl}[1]{\mbox{\textup{\boldmath$#1$}}}
\renewcommand{\Bigr}[1]{\mbox{\textup{\boldmath$#1$}}}

\newcommand{\m}{\mathrm{m}}
\newcommand{\p}{\mathrm{p}}
\newcommand{\Spec}{\mathop{\mathrm{Spec}}\nolimits}
\newcommand{\Max}{\mathop{\mathrm{Max}}\nolimits}
\newtheorem{remarkaux}[theorem]{{\textit{Remark}}}
\newenvironment{remark}{\begin{remarkaux}\upshape}{\end{remarkaux}}

\begin{document}

\begin{center}
\textbf{\LARGE
     Elementary Factors and Reduced Minors\\
     for Linear Systems over Commutative Rings
}\medskip\\
{\large               Kazuyoshi~MORI$\dag$
}
\\
       Department of Electrical Engineering, Tohoku University\\
       Aoba-ku, Aramaki aza-Aoba 05, Sendai 980-8579, JAPAN\\
         (\texttt{Kazuyoshi.MORI@IEEE.ORG})
\end{center}

\renewcommand{\thefootnote}{\fnsymbol{footnote}}
\footnotetext[2]{
  This work was partially done while the 
  author was visiting the
  Institut de Recherche en Cybern\'etique de Nantes, France and the
  author thanks this institution for its hospitality and
  support.}
\renewcommand{\thefootnote}{\arabic{footnote}}

  \centerline{\textsl{March 10, 2000}}

\section*{Abstract}
In 1994, Sule presented the necessary and sufficient conditions of the
feedback stabilizability of systems over unique factorization domains
in terms of elementary factors and in terms of reduced minors.
Recently, Mori and Abe have generalized his theory over commutative
rings.  They have introduced the notion of the generalized elementary
factor, which is a~generalization of the elementary factor, and have
given the necessary and sufficient condition of the feedback
stabilizability.  In this paper, we present two generalization of the
reduced minors.  Using each of them, we state the necessary and
sufficient condition of the feedback stabilizability over commutative
rings.  Further we present the relationship between the
generalizations and the generalized elementary factors.

\section*{Keywords}
  Linear systems,  Feedback stabilization, Factorization approach,
  Systems over rings

\section{Introduction}\label{S:Introduction}
This paper is concerned with the coordinate-free approach to control
systems.  The coordinate-free approach is a~factorization approach but
does not require the coprime factorizations of plants.

The factorization approach was patterned after Desoer \emph{et
al.}\cite{bib:desoer80a} and Vidyasagar \emph{et
al.}\cite{bib:vidyasagar82a}, which has the advantage that it
embraces, within a~single framework, numerous linear systems such as
continuous-time as well as discrete-time systems, lumped as well as
distributed systems, $1$-D as well as~$n$-D (multidimensional)
systems, etc.\cite{bib:vidyasagar82a}.  In this approach, when
problems such as feedback stabilization are studied, one can {focus}
on the key aspects of the problem under study rather than be
distracted by the special features of a~particular class of linear
systems.  A~transfer function of this approach is considered as the
ratio of two stable causal transfer functions and the set of stable
causal transfer functions forms a~commutative ring.  For a~long time,
the theory of the factorization approach had been founded on the
coprime factorizability of transfer matrices, which is satisfied in
the case where the set of stable causal transfer functions is such
a~commutative ring as a~Euclidean domain, a~principal ideal, or
a~B\'ezout domain.

However, Anantharam in~\cite{bib:anantharam85a} showed that there
exist models in which some stabilizable plants do not have
right-/left-coprime factorizations.  He considered the case where
$\Zset[\sqrt{5}\mathrm{i}]$ ($\simeq\Zset[x]/(x^2+5)$) is the set of
stable causal transfer functions, where $\Zset$ is the ring of
integers and~$\mathrm{i}$ the imaginary unit.  Using it, he showed
that there exists a~stabilizable plant which does not have
right-/left-coprime factorizations.  Further Mori
in~\cite{bib:mori99b} has recently considered the case where
$\Rset[z^2,z^3]$ is the set of stable causal transfer functions,
where~$z$ denotes the unit delay operator and $\Rset$ the real field.
This set is corresponding to the discrete finite-time delay system
which does not have the unit delay.  He has presented that in the
model, some stabilizable plants do not have right-/left-coprime
factorizations.  Both $\Zset[\sqrt{5}\mathrm{i}]$ and $\Rset[z^2,z^3]$
are not unique factorization domains.

Sule in~\cite{bib:sule94a,bib:sule98a} has presented a~theory of the
feedback stabilization of multi-input multi-output strictly causal
plants over commutative rings with some restrictions.  This approach
to the stabilization theory is called ``coordinate-free approach'' in
the sense that the coprime factorizability of transfer matrices is not
required.  

In the case where the set of stable causal transfer functions is
a~unique factorization domain, Sule in~\cite{bib:sule94a} introduced
two notions, that is, elementary factors and reduced minors.  Using
each of them he gave the necessary and sufficient condition of the
feedback stabilizability of the causal plants over commutative rings
(Theorem\,4 and Corollary\,2 of~\cite{bib:sule94a}).  Especially,
using elementary factors, Sule presented a~construction method of
a~stabilizing controller of a~stabilizable plant.  Recently, Mori and
Abe in~\cite{bib:mori97a,bib:mori98bsnogzip} have generalized his
theory over commutative rings.  They have introduced the notion of the
generalized elementary factor, which is a~generalization of the
elementary factor, and have given the necessary and sufficient
condition of the feedback stabilizability.  Further Lin
in~\cite{bib:lin98b} has presented the necessary and sufficient
condition of the (structural) stabilizability of the multidimensional
systems with the construction method of a~stabilizing controller.  In
the case of the structural stability\cite{bib:guiver85a}, it is known
that the set of stable causal transfer functions is a~unique
factorization domain.  Lin in~\cite{bib:lin98b} introduced a~notion
``generating polynomial'' about the plants and presented the necessary
and sufficient condition of the stabilizability of the
multidimensional systems with the construction method of a~stabilizing
controller.  It is known that the notion of the generating polynomial
is equivalent to the notion of the reduced minors.

In this paper we have two main objectives.  The first one is to
generalize the notion of the reduced minors and, using the
generalizations, to state the necessary and sufficient condition of
the feedback stabilizability over commutative rings since the original
definition has been given on unique factorization domains.  We will
present two generalizations.  The other is to present the relationship
between the generalizations and the generalized elementary factors.

Historically the minors concerning the plants are much
investigated (e.g.~%
  \cite{bib:berenstein86a,
        bib:lin88a,
        bib:lin88b,
        bib:lin92a,
        bib:lin99a,
        bib:wood98a,
        bib:youla79a,
        bib:youla84a,
        bib:zerz96a}).
We will present that in the coordinate-free approach, the minors can
play a~role to state the feedback stabilizability, that is, the
\emph{projectivity} of the ideal generated by minors concerning the
plant is a~criterion of the feedback stability.

This paper is organized as follows.  After this introduction, we begin
on the preliminary in Section\,\ref{S:Preliminary}, in which we give
mathematical preliminaries, set up the feedback stabilization problem
and present the previous results.  In Section\,\ref{S:PR}, we present
the previous results of the feedback stabilizability expressed with
the elementary factors, its derivation, and the reduced minors.  We
present a generalization of the reduced minor in Section\,\ref{S:4}
and using it present the necessary and sufficient condition of the
feedback stabilizability over commutative rings in Section\,\ref{S:5}.
Then in Section\,\ref{S:6} we present another generalization of the
reduced minors and its relation to the generalized elementary factors.

\section{Preliminaries}
\label{S:Preliminary}
In the following we begin by introducing the notations of commutative
rings, matrices, and modules used in this paper.  Then we give the
formulation of the feedback stabilization problem.

\subsection{Notations}\label{SS:Notations}
\paragraph{Commutative Rings}
In this paper, we consider that any commutative ring has the
identity~$1$ different from zero.  Let~${\cal R}$ denote a~(unspecified)
commutative ring.  The total ring of fractions of~${\cal R}$ is denoted by~$
{\cal F}({\cal R})$.

We will consider that \emph{the set of stable causal transfer
functions} is a~commutative ring, which is denoted by~${\cal A}$ throughout
this paper.  Further, we will use the following rings of fractions.
\begin{romannum}
\item
The first one appears as the total ring of fractions of~${\cal A}$, which is
denoted by~${\cal F}({\cal A})$ or simply by~${\cal F}$; that is, ${\cal F}=\{ n/d\,|\, n,d\in{\cal A},~
\mbox{$d$ is a~nonzerodivisor}\}$.  This will be considered as
\emph{the set of all possible transfer functions}.
\item
Let~$f$ denote a~nonzero (but possibly nonzerodivisor) element of~${\cal A}
$.  Given a~set $S_f=\{ 1,f,f^2,\ldots\}$, which is a~multiplicative
subset of~${\cal A}$, we denote by~${\cal A}_f$ the ring of fractions of~${\cal A}$ with
respect to the multiplicative subset~$S_f$; that is, ${\cal A}_f=\{ n/d\,|\, n\in
{\cal A},~d\in S_f\}$.
\item
Let~$\p$ denote a~prime ideal of~${\cal A}$ and~$S$ the
complement of the prime ideal~$\p$, that is, $S={\cal A}\backslash\p$.  Then~$S$ is
a~multiplicative subset of~${\cal A}$.  We denote by~${\cal A}_\p$ the ring of
fractions of~${\cal A}$ with respect to the multiplicative subset~$S$; that
is, ${\cal A}_\p =\{ n/d\,|\, n\in{\cal A},~d\in S\}$.
\item
The last one is the total ring of fractions of~${\cal A}_f$ or ${\cal A}_\p$,
which is denoted by ${\cal F}({\cal A}_f)$ and ${\cal F}({\cal A}_\p)$; that is,
  ${\cal F}({\cal A}_f)=
   \{ n/d\,|\, n,d\in{\cal A}_f,~\mbox{$d$ is a~nonzerodivisor of ${\cal A}_f$}\}$
and
  ${\cal F}({\cal A}_\p)=
   \{ n/d\,|\, n,d\in{\cal A}_\p,~\mbox{$d$ is a~nonzerodivisor of ${\cal A}_\p$}\}$.
If~$f$ is a~nonzerodivisor of~${\cal A}$, ${\cal F}({\cal A}_f) $ coincides with the
total ring of fractions of~${\cal A}$.  Otherwise, they do not coincide.
\end{romannum}

In the case where~${\cal A}$ is a~unique factorization domain, we call~$a$
in~${\cal A}$ \emph{the radical of} $b$ in~${\cal A}$ if~$a$ has all nonunit
factors of~$b$ and is squarefree, that is, $a$ does not have
duplicated nonunit factors.  Note here that the radical defined here
is unique up to any unit multiple.

For convenience, throughout the paper, if $a\in{\cal A}$ ($a\in{\cal R}$), then~$a$
itself denotes $a/1$ in~${\cal A}_f$ and ${\cal A}_\p$ ($a/1$ in ${\cal F}({\cal R})$).  Moreover
if $a\in{\cal A}_f$ or ${\cal A}_\p$ ($a\in{\cal R}$) and if there exists $b\in{\cal A}$ such that
$a=b/1$ over~${\cal A}_f$ or ${\cal A}_\p$ (over ${\cal F}({\cal R})$), then we regard~$a$ as an
element of~${\cal A}$ (${\cal R}$).

In the rest of the paper, we will use~${\cal R}$ as an unspecified
commutative ring and mainly suppose that~${\cal R}$ denotes one of ~${\cal A}$, $
{\cal A}_f$, and ${\cal A}_\p$.

We will denote by~$\Spec({\cal R})$ the set of all prime ideals of~${\cal R}$ and
by~$\Max({\cal R})$ the set of all maximal ideals of~${\cal R}$.  Suppose that
$\mfa$ is an ideal of~${\cal R}$.  Then we denote by $\mfa_f$ the ideal of
fractions of $\mfa$ with respect to
  $\{ 1,f,f^2,\ldots\}$ with $f\in{\cal R}$ 
       (that is, $\mfa_f=\{ n/d\,|\, n\in\mfa,~d\in\{ 1,f,f^2,\ldots\}\}$)
and by $\mfa_\p$ the ideal of fractions of $\mfa$ with respect to ${\cal R}\backslash
\p$ with $\p\in\Spec({\cal R})$ (that is, $\mfa_\p=\{ n/d\,|\, n\in\mfa,~d\in{\cal R}\backslash\p\}
$).  If $\mfa$ is an ideal of~${\cal R}$ and if~$S$ is a~subset of~${\cal R}$,
then we denote by $(\mfa:S)$ the \emph{quotient ideal} which is the
set $\{ f\in{\cal R}\,|\, fS\subset\mfa\}$.

The reader is referred to~Chapter\,3 of~\cite{bib:atiyah69a} for the
ring of fractions.

\paragraph{Matrices}
The set of matrices over~${\cal R}$ of size $x\times y$ is denoted by ${\cal R}^{x\times
y}$.  Further, the set of square matrices over~${\cal R}$ of size~$x$ is
denoted by~$ ({\cal R})_x$.  The identity and the zero matrices are denoted
by~$E_x$ and $O_{x\times y}$, respectively, if the sizes are required,
otherwise they are denoted by~$E$ and~$O$.

Matrix~$A$ over~${\cal R}$ is said to be \emph{nonsingular}
$\Bigl($\emph{singular}$\Bigr)$ \emph{over~${\cal R}$} if the determinant of
the matrix~$A$ is a~nonzerodivisor $\Bigl($a zerodivisor$\Bigr)$ of~$
{\cal R}$.  Matrices~$A$ and~$B$ over~${\cal R}$ are \emph{right-}
\emph{$\Bigl($left-$\Bigr)$coprime over~${\cal R}$} if there exist
matrices~$X$ and~$Y$ over~${\cal R}$ such that
  $XA+YB=E$ $\Bigl(AX+BY=E\Bigr)$ 
holds.  Note that, in the sense of the above definition, two matrices
which have no common right-$\Bigl($left-$\Bigr)$factors except
invertible matrices may not be right-$\Bigl($left-$\Bigr)$coprime
over~${\cal R}$.  Further, an ordered~pair $(N,D)$ of matrices~$N$ and~$D$
is said to be a \emph{right-coprime factorization over~${\cal R}$} of~$P$ if
(i)~$D$ is nonsingular over~${\cal R}$, (ii) $P=ND^{-1}$ over ${\cal F}({\cal R})$, and
(iii)~$N$ and~$D$ are right-coprime over~${\cal R}$.  As the parallel
notion, the \emph{left-coprime factorization over~${\cal R}$} of~$P$ is
defined analogously.  That is, an ordered~pair
$(\widetilde{D},\widetilde{N})$ of matrices~$\widetilde{N}$
and~$\widetilde{D}$ is said to be a \emph{left-coprime factorization
over~${\cal R}$} of~$P$ if (i) $\widetilde{D}$ is nonsingular over~${\cal R}$,
(ii) $P=\widetilde{D}^{-1}\widetilde{N}$ over ${\cal F}({\cal R})$, and
(iii)~$\widetilde{N}$ and~$\widetilde{D}$ are left-coprime over~${\cal R}$.
Note that the order of the ``denominator'' and ``numerator'' matrices
is interchanged in the latter case.  This is to reinforce the point
that if $(N,D)$ is a~right-coprime factorization over~${\cal R}$ of~$P$,
then $P=ND^{-1}$, whereas if $(\widetilde{D},\widetilde{N})$ is
a~left-coprime factorization over~${\cal R}$ of~$P$, then
$P=\widetilde{D}^{-1}\widetilde{N}$ according
to~\cite{bib:vidyasagar85a}.  For short, we may omit ``over~${\cal R}$''
when ${\cal R}={\cal A}$, and ``right'' and ``left'' when the size of matrix is $1
\times 1$.  In the case where matrices are potentially used to express
\emph{left} fractional form and/or \emph{left} coprimeness, we usually
attach a~tilde `$\widetilde{~~}$' to symbols; for example
$\widetilde{N}$, $\widetilde{D}$ for
$P=\widetilde{D}^{-1}\widetilde{N}$ and $\widetilde{Y}$,
$\widetilde{X}$ for $\widetilde{Y}N+\widetilde{X}D=E$.

\paragraph{Modules}
Let $M_r(X)$ \,$\Bigl(M_c(X)\Bigr)$ denote the~${\cal R}$-module generated
by rows \,$\Bigl($columns$\Bigr)$ of a~matrix~$X$ over~${\cal R}$.
Let $X=AB^{-1}=\widetilde{B}^{-1}\widetilde{A}$ be a~matrix over~${\cal F}
({\cal R})$, where~$A$,~$B$,~$\widetilde{A}$,~$\widetilde{B}$ are matrices
over~${\cal R}$.  It is known that $M_r(\bmatrix{A^t & B^t}^t)$
\,\,$\Bigl(M_c(\bmatrix{\widetilde{A} & \widetilde{B}})\Bigr)$ is
unique up to an isomorphism with respect to any choice of fractions
$AB^{-1}$ of~$X$ \,$\Bigl(\widetilde{B}^{-1}\widetilde{A}$ of
$X\Bigr)$ (Lemma\,2.1 of\,\cite{bib:mori97a}).  Therefore, for
a~matrix~$X$ over~${\cal R}$, we denote by~${\cal T}_{X,{\cal R}}$ and~${\cal W}_{X,{\cal R}}$ the
modules $M_r(\bmatrix{A^t & B^t}^t)$ and $M_c(\bmatrix{\widetilde{A} &
\widetilde{B}})$, respectively.

An~${\cal R}$-module~$M$ is called \emph{free} if it has a~basis, that is,
a~linearly independent system of generators.  The \emph{rank} of
a~free~${\cal R}$-module~$M$ is equal to the cardinality of a~basis of~$M$,
which is independent of the basis chosen.  An~${\cal R}$-module~$M$ is
called \emph{projective} if it is a~direct summand of a~free~${\cal R}
$-module, that is, there is a~module~$N$ such that $M\oplus N$ is free.
The reader is referred to~Chapter\,2 of~\cite{bib:atiyah69a} for the
module theory.

We will consider occasionally ideals as modules in this paper.  So, we
will apply the words ``projective,'' ``free,'' and ``isomorphic'' to
ideals.  It is easy to check that an ideal which is free as a~module
is equivalent to a~principal ideal whose generator is
a~nonzerodivisor.

\subsection{Feedback Stabilization Problem}\label{SS:FSProblem}
\begin{figure}
\begin{center}
\setlength{\unitlength}{0.07mm}
\begin{picture}(1160.000000,400.000000)(-50.000000,-50.000000)
\linethickness{0.4mm}
\put(220.000000,50.000000){\framebox(200.000000,150.000000){\LARGE$C$}}
\put(640.000000,50.000000){\framebox(200.000000,150.000000){\LARGE$P$}}
\put(110.000000,125.000000){\circle{30.000000}}
\put(530.000000,125.000000){\circle{30.000000}}
\put(0.000000,125.000000){\vector(1,0){95.000000}}
\put(95.000000,128.000000){\vector(1,0){0}}
\put(95.000000,122.000000){\vector(1,0){0}}
\put(125.000000,125.000000){\vector(1,0){95.000000}}
\put(220.000000,128.000000){\vector(1,0){0}}
\put(220.000000,122.000000){\vector(1,0){0}}
\put(420.000000,125.000000){\vector(1,0){95.000000}}
\put(515.000000,128.000000){\vector(1,0){0}}
\put(515.000000,122.000000){\vector(1,0){0}}
\put(545.000000,125.000000){\vector(1,0){95.000000}}
\put(640.000000,128.000000){\vector(1,0){0}}
\put(640.000000,122.000000){\vector(1,0){0}}
\put(840.000000,125.000000){\vector(1,0){190.000000}}
\put(1030.000000,128.000000){\vector(1,0){0}}
\put(1030.000000,122.000000){\vector(1,0){0}}
\put(110.000000,0.000000){\vector(0,1){110.000000}}
\put(107.000000,110.000000){\vector(0,1){0}}
\put(113.000000,110.000000){\vector(0,1){0}}
\put(110.000000,0.000000){\line(1,0){810.000000}}
\put(920.000000,0.000000){\line(0,1){125.000000}}
\put(530.000000,300.000000){\vector(0,-1){160.000000}}
\put(527.000000,140.000000){\vector(0,-1){0}}
\put(533.000000,140.000000){\vector(0,-1){0}}
\put(0.000000,135.000000){\makebox(0,0)[lb]{\Large$u_1$}}
\put(218.000000,140.000000){\makebox(0,0)[rb]{\Large$e_1$}}
\put(90.000000,115.000000){\makebox(0,0)[rt]{\Large$+$}}
\put(120.000000,85.000000){\makebox(0,0)[lt]{\Large$-$}}
\put(430.000000,135.000000){\makebox(0,0)[lb]{\Large$y_1$}}
\put(520.000000,290.000000){\makebox(0,0)[rt]{\Large$u_2$}}
\put(638.000000,140.000000){\makebox(0,0)[rb]{\Large$e_2$}}
\put(1028.000000,135.000000){\makebox(0,0)[rb]{\Large$y_2$}}
\put(510.000000,115.000000){\makebox(0,0)[rt]{\Large$+$}}
\put(540.000000,205.000000){\makebox(0,0)[lb]{\Large$+$}}
\put(920.000000,125.000000){\circle*{22.500000}}
\end{picture}

  \caption{Feedback system~$\Sigma$.}\label{Fig:FeedbackSystem}
\end{center}
\end{figure}
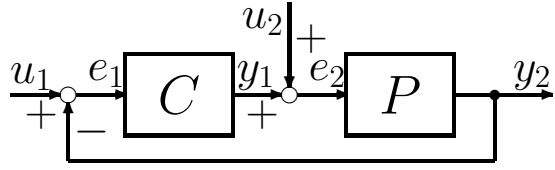

The stabilization problem considered in this paper follows that of
Sule in~\cite{bib:sule94a}, and Mori and Abe
in~\cite{bib:mori97a}, who consider the feedback system~$
\Sigma$~\cite[Ch.5, Figure\,5.1]{bib:vidyasagar85a} as in
Figure\,\ref{Fig:FeedbackSystem}.
For further details the reader is referred
to~\cite{bib:vidyasagar85a}.  Throughout the paper, the plant we
consider has~$m$ inputs and~$n$ outputs, and its transfer matrix,
which is also called a \emph{plant} itself simply, is denoted by~$P$
and belongs to~${\cal F}^{n\times m}$.  We can always represent~$P$ in the form
of a~fraction $P=ND^{-1}$
$\Bigl(P=\widetilde{D}^{-1}\widetilde{N}\Bigr)$, where $N\in{\cal A}^{n\times m}$
\,$\Bigl(\widetilde{N}\in{\cal A}^{n\times m}\Bigr)$ and $D\in({\cal A})_m$
\,$\Bigl(\widetilde{D}\in({\cal A})_n\Bigr)$ with nonsingular~$D$
\,$\Bigl(\widetilde{D}\Bigr)$.  
\begin{definition}
\label{D:Stability}
For $P\in{\cal F}^{n\times m}$ and $C\in{\cal F}^{m\times n}$, a~matrix $H(P,C)\in({\cal F})_{m+n}$ is
defined as
\begin{equation}\label{E:H(P,C)}
 ~~~~~~~~
  H(P,C) \!=\!
   \bmatrix{
    (E_n+PC)^{-1}  &  -P(E_m+CP)^{-1} \cr
    C(E_n+PC)^{-1} & (E_m+CP)^{-1}}
\end{equation}
provided that $\det(E_n+PC)$ is a~nonzerodivisor of~${\cal A}$.  This
$H(P,C)$ is the transfer matrix from $\bmatrix{u_1^t & u_2^t}^t$ to
$\bmatrix{e_1^t & e_2^t}^t$ of the feedback system~$\Sigma$.  If (i)
$\det(E_n+PC)$ is a~nonzerodivisor of~${\cal A}$ and (ii) $H(P,C)\in({\cal A})_
{m+n}$, then we say that the plant~$P$ is \emph{stabilizable}, $P$ is
\emph{stabilized} by~$C$, and~$C$ is a \emph{stabilizing controller}
of~$P$.
\end{definition}

Since the transfer matrix $H(P,C)$ of the stable causal feedback
system has all entries in~${\cal A}$, we call the above notion~\emph{${\cal A}
$-stabilizability}.  One can further introduce the notion of~\emph{${\cal R}
$-stabilizability} with either ${\cal R}={\cal A}_f$ or $={\cal A}_\p$ as follows.
\begin{definition}
\label{D:RStability}
Suppose that~${\cal R}$ is either~${\cal A}_f$ with $f\in{\cal A}\backslash\{ 0\}$ or ${\cal A}_\p$ with
$\p\in\Spec({\cal A})$.  If (i) $\det(E_n+PC)$ is a~nonzerodivisor of~${\cal R}$
and (ii) $H(P,C)\in({\cal R})_ {m+n}$, then we say that the plant~$P$ is
\emph{${\cal R}$-stabilizable}, $P$ is \emph{${\cal R}$-stabilized} by~$C$,
and~$C$ is an \emph{${\cal R}$-stabilizing controller} of~$P$.
\end{definition}

The causality of transfer functions is an important physical
constraint.
We employ, in this paper, the definition of the causality from
Vidyasagar~\emph{et al.}\cite[Definition\,3.1]{bib:vidyasagar82a}.
\begin{definition}\label{D:Causality}
Let~${\cal Z}$ be a~prime ideal of~${\cal A}$, with ${\cal Z}\neq{\cal A}$, including all
zerodivisors.  Define the subsets~${\cal P}$ and ${\cal P}_{\textrm{s}}$ of~${\cal F}$
as follows:
\begin{eqnarray*}
  {\cal P}&=&\{ a/b\in{\cal F}\,|\, a\in{\cal A},~ b\in{\cal A}\backslash{\cal Z}\},\\
  {\cal P}_{\textrm{s}}&=&\{ a/b\in{\cal F}\,|\, a\in{\cal Z},~ b\in{\cal A}\backslash{\cal Z}\}.
\end{eqnarray*}
Then every transfer function in~${\cal P}$ $\Bigl({\cal P}_{\textrm{s}}\Bigr)$ is
called \emph{causal} $\Bigl($\emph{strictly causal}$\Bigr)$.
Analogously, if every entry of a~transfer matrix~$F$ is in~${\cal P}$
$\Bigl({\cal P}_ {\textrm{s}}\Bigr)$, the transfer matrix~$F$ is called
\emph{causal} $\Bigl($\emph{strictly causal}$\Bigr)$.  A~matrix over~$
{\cal A}$ is said to be \emph{${\cal Z}$-nonsingular} if the determinant is in ${\cal A}\backslash
{\cal Z}$, and \emph{${\cal Z}$-singular} otherwise.
\end{definition}

Before proceeding the next section, we here introduce several symbols
used throughout this paper.  The symbol~${\cal I}$
denotes the family of all sets of~$m$ distinct integers between~$1$
and $m+n$, and~${\cal J}$
the family of all sets of~$n$ distinct integers between~$1$ and $m+n$
(recall that~$m$ and~$n$ are the numbers of the inputs and the
outputs, respectively).  Normally, elements of~${\cal I}$ (${\cal J}$) will be
denoted by~$I$ ($J$) possibly with suffices.  They will be used as
suffices as well as sets.  If~$I$ is an element of~${\cal I}$ and if
$i_1,\ldots, i_m$ are elements of~$I$ with ascending order, that is,
$i_a<i_b$ if $a<b$, then the symbol~$\Delta_I$ denotes the $m\times(m+n)$
matrix whose $(k,i_k)$-entry is~$1$ for $i_k\in I$ and zero otherwise.
Analogously if~$J$ is an element of~${\cal J}$ and if $j_1,\ldots, j_n$ are
elements of~$J$ with ascending order, then the symbol~$\Delta_J$ denotes
the $n\times(m+n)$ matrix whose $(k,j_k)$-entry is~$1$ for $j_k\in J$ and zero
otherwise.

\section{Previous Results}
\label{S:PR}
In this section, we recall the previous results about the necessary
and sufficient condition of the feedback stabilizability.  First one
is stated in terms of the elementary factors and the other in terms of
the reduced minors.

\subsection{Feedback Stabilizability in terms of Elementary Factors}
\label{SS:20.Jun.99.190807}
To state the result, we first recall the notion of the elementary
factors, which was defined under the assumption that~${\cal A}$ is a~unique
factorization domain.

\begin{definition}\label{D:EF}
\textup{(Elementary Factors, \cite[p.1689]{bib:sule94a})~} Suppose
that~${\cal A}$ is a~unique factorization domain.  Denote by~$T$ and~$W$ the
matrices $\bmatrix{N^t&dE_m}^t$ and $\bmatrix{N&dE_n}^t$ over~${\cal A}$
with $P=Nd^{-1}$.  Further denote by ${\cal I}^*$ $\Bigl({\cal J}^*\Bigr)$ the set
of $I$'s in~${\cal I}$ $\Bigl(J$'s in ${\cal J}\Bigr)$ such that $\Delta_IT$ $\Bigl(\Delta_
JW^t\Bigr)$ is nonsingular.  Then for each $I\in{\cal I}^*$, let~$f_I$ be the
radical of the least common multiple of all the denominators of the
matrix $T (\Delta_ IT)^{-1}$ and for each $J\in{\cal J}^*$, $g_J$ be the radical
of the least common multiple of all the denominators of the matrix
$W^t (\Delta_ JW^t)^{-1}$.  Then~$f_I$~($g_J$) is called \emph{the
elementary factor of the matrix~$T$ $\Bigl(W\Bigr)$ with respect to $I
\in{\cal I}$} $\Bigl(J\in{\cal J}\Bigr)$, $F=\{ f_I\,|\, I\in{\cal I}^*\}$ \emph{the family of
elementary factors of the matrix $T$}, $G=\{ g_J\,|\, J\in{\cal J}^*\}$ \emph{the
family of elementary factors of the matrix $W$}, and $H=
\{ h_{IJ}:=f_Ig_J\,|\, I\in{\cal I}^*, J\in{\cal J}^*\}$ \emph{ the family of elementary
factors of $P$}.
\end{definition}

Then the necessary and sufficient condition of the feedback
stabilizability is given as follows.
\begin{theorem}\label{Th:21.Jun.99.190413}
\textup{\rm (Theorem\,4 of~\cite{bib:sule94a})~} Suppose that~${\cal A}$ is
a~unique factorization domain.  Then the plant~$P$ is stabilizable if
and only if the elementary factors of~$P$ are coprime, that is,
  $\sum_{I\in{\cal I}^*, J\in{\cal J}^*}(h_{IJ})={\cal A}$. 
\end{theorem}

In the proof of this theorem, Sule gave a~method to construct
a~stabilizing controller of the plant.

The result above has been extended to include systems over commutative
rings by Mori and Abe in~\cite{bib:mori98bsnogzip} as follows.  They
introduced the notion of the generalized elementary factors, which is
a~generalization of the elementary factors, and using it, stated the
necessary and sufficient conditions of the feedback stabilizability
over commutative rings.

\begin{definition}\label{D:GEF}
\textup{(Generalized Elementary Factors, Definition\,3.1
of~\cite{bib:mori98bsnogzip})~} Denote by~$T$ the matrix
$\bmatrix{N^t&D^t}^t$ over~${\cal A}$ with $P=ND^{-1}$.  For each $I\in{\cal I}$,
an ideal $\Lambda_{P\!I}$ over~${\cal A}$ is defined as
\[
  \Lambda_{P\!I}=\{\lambda\in{\cal A}\,|\,\exists K\in{\cal A}^{(m+n)\times m}\  \lambda T  =K \Delta_I T\}.
\]
We call the ideal $\Lambda_{P\!I}$ the \emph{generalized elementary factor} of
the plant~$P$ with respect to~$I$.  Further, the set of all $\Lambda_{P\!I}$'s
is denoted by~${\cal L}_{P}$, that is, ${\cal L}_{P}=\{\Lambda_{P\!I}\,|\, I\in{\cal I}\}$.
\end{definition}

In the case where~${\cal A}$ is a~unique factorization domain, a~generalized
elementary factor with respect to $I\in{\cal I}$ is a~principal ideal and the
radical of its generator is an~elementary factor of~$T$ with respect
to~$I$ up to a~unit multiple.

It is known that the generalized elementary factor of a~plant~$P$ is
independent of the choice of fractions $ND^{-1}=P$ (Lemma\,3.3 of
\cite{bib:mori98bsnogzip}).

The following is the necessary and sufficient conditions of the
feedback stabilizability.
\begin{theorem}\label{Th:3.5}
\textup{\rm (Theorem\,3.2 of~\cite{bib:mori98bsnogzip})~} Consider
a~causal plant~$P$.  Then the following statements are equivalent:
\begin{romannum}
\item
The plant~$P$ is stabilizable.
\item
${\cal A}$-modules~${\cal T}_{P,{\cal A}}$ and~${\cal W}_{P,{\cal A}}$ are projective.
\item
The set of all generalized elementary factors of~$P$ generates~${\cal A}$;
that is, ${\cal L}_P$ satisfies:
\begin{equation}\label{E:Th:3.5:1}
  \sum_{\Lambda_{P\!I}\in{\cal L}_P} \Lambda_{P\!I}={\cal A}.
\end{equation}
\end{romannum}
\end{theorem}

Provided that we can check (\ref{E:Th:3.5:1}) and that we can
construct the right-coprime factorizations over ${\cal A}_{\lambda_I}$ of the
given causal plant, where~$\lambda_I$ is a~nonzero element of~${\cal A}$, Mori
and Abe\cite{bib:mori98bsnogzip} have given a~method to construct
a~causal stabilizing controller of a~causal stabilizable plant, which
has been given in the proof of ``(iii)$\rightarrow$(i)'' of Theorem\,3.2
of~\cite{bib:mori98bsnogzip}.

\subsection{Feedback Stabilizability in terms of Reduced Minors}
\label{SS:21.Jun.99.183641}
We first recall the definition of the reduced minors and then state
the necessary and sufficient conditions of the feedback
stabilizability in terms of the reduced minors.  We suppose in this
subsection that~${\cal A}$ is a~unique factorization domain.

\begin{definition}\label{D:ReducedMinors}
\textup{\rm (Reduced Minors, \cite[p.1690]{bib:sule94a})~} Let~$P$ be
a~plant of ${\cal F}^{n\times m}$, $N$ a~matrix of ${\cal A}^{n\times m}$, and~$d$ an
element of~${\cal A}$ such that $P=Nd^{-1}$.  Denote by~$T$ and~$W$ the
matrices $\bmatrix{N^t&dE_m}^t$ and $\bmatrix{N&dE_n}$.  Let
$t_I=\det(\Delta_IT)$ $\Bigl(w_J=\det(\Delta_JW^t)\Bigr)$, which is
a~full-size minor of the matrix~$T$ $\Bigl(W\Bigr)$, for $I\in{\cal I}$
$\Bigl(J\in{\cal J}\Bigr)$.  Let~$d_t$ $\Bigl(d_w\Bigr)$ be the greatest
common factor of~$t_I$'s $\Bigl(w_J$'s$\Bigr)$ and $a_I=t_I/d_t$ for
$I\in{\cal I}$ $\Bigl(b_J=w_J/d_w$ for $J\in{\cal J}\Bigr)$.  Then~$a_I$
$\Bigl(b_J\Bigr)$ is called the \emph{reduced minor of the matrix~$T$
($W$) with respect to $I\in{\cal I}$ $\Bigl(J\in{\cal J}\Bigr)$}, the set $\{ a_I\,|\, I
\in{\cal I}\}$ $\Bigl(\{ b_J\,|\, J\in{\cal J}\}\Bigr)$ \emph{the family of reduced minors
of~$T$ $\Bigl(W\!\Bigr)$}.
\end{definition}

It is known that the families of reduced minors of~$T$ and of~$W$ are
identical modulo units (Lemma\,5 of~\cite{bib:sule94a}).

Now, Corollary\,2 of~\cite{bib:sule94a} including its comments can be
stated as follows:
\begin{theorem}\label{Th:3.4}
\textup{\rm (cf.~Corollary\,2 of~\cite{bib:sule94a})~} Suppose that~$
{\cal A}$ is a~unique factorization domain.  A~plant $P\in{\cal F}^{m\times n}$ is
stabilizable if and only if the family of the reduced minors of~$T$
(and also of $W$) generates~${\cal A}$.
\end{theorem}

The theorem above can be rewritten directly as follows.
\begin{corollary}
\label{R:19.Jul.99.204321}
Let~$t_I$ and~$w_J$ be as in
Definition\ref{D:ReducedMinors}.
Then the following are equivalent:
\begin{romannum}
\item A~plant $P\in{\cal F}^{m\times n}$ is stabilizable.
\item The ideal $\sum_{I\in{\cal I}}(t_I)$ is principal, or
equivalently free as an ${\cal A}$-module.
\item The ideal $\sum_{J\in{\cal J}}(w_J)$ is principal, or
equivalently free as an ${\cal A}$-module.
\end{romannum}
\end{corollary}

\section{Full-Size Minor Ideal}\label{S:4}
On the statements concerning the elementary factors and the reduced
minors in Subsections\,\ref{SS:20.Jun.99.190807}
and~\ref{SS:21.Jun.99.183641}, we have considered that the denominator
matrices of the plant is expressed as $dE_m$ or $dE_n$ rather than
general nonsingular matrices.  This may be considered as a~restriction
on the expression of the plant.  Thus we rather consider that~$P$ is
expressed as either $P=ND^{-1}$ with $N\in{\cal A}^{n\times m}$ and $D\in({\cal A})_m$ or
$P=\widetilde{D}^{-1}\widetilde{N}$ with $\widetilde{N}\in{\cal A}^{n\times m}$
and $\widetilde{D}\in({\cal A})_n$.  Now we redefine the matrices~$T$, $W$ as
$T=\bmatrix{ N^t&D^t}^t$ and $W=\bmatrix{
\widetilde{N}&\widetilde{D}}$.  Further we consider that $t_I$'s and
$w_J$'s are defined with the matrices~$T$ and~$W$ here.  In the rest
of this paper, we will use these notations unless otherwise stated.

We now introduce a~notion to state the feedback stabilizability over
commutative rings.
\begin{definition}\label{D:FullSizeMinorIdeals}
\textup{(Full-Size Minor Ideals)~}
The ideal generated by $t_I$'s for $I\in{\cal I}$ is called the
\emph{full-size minor ideal} of the plant~$P$.  We denote it by $\sum_
{I\in{\cal I}}(t_I)$ or simply $\mft$.
\end{definition}

We can also consider the ideal generated by $w_J$'s for $J\in{\cal J}$,
denoted by $\sum_{J\in{\cal J}}(w_J)$ or simply~$\mfw$.  The ideals $\mft$
and $\mfw$ depend on the fractional representation of the plant
$P=ND^{-1}=\widetilde{D}^{-1}\widetilde{N}$.  However, this is
\emph{not} a~problem from the following reason.  To state the feedback
stabilizability in terms of the full-size minor ideals, we will regard
them as modules.  Further, when these ideals are considered as
modules, both the ideals $\mft$ and $\mfw$ are uniquely determined as
modules up to isomorphism with respect to any choice of fractions
$ND^{-1}$ and $\widetilde{N}^{-1}\widetilde{D}$ of~$P$ as shown below.

\begin{lemma}\label{L:4.2}
Let~$P$ be in ${\cal F}({\cal R})^{n\times m}$, where~${\cal R}$ is one of~${\cal A}$, ${\cal A}_f$ with
a~nonzero~$f\in{\cal A}$, and ${\cal A}_\p$ with a~prime ideal~$\p$ in $\Spec({\cal A})$.
For $x=1,2$ let
  $N_x$, $D_x$,
  $\widetilde{N}_x$, $\widetilde{D}_x$ be 
matrices over~${\cal R}$ with
$P=N_xD_x^{-1}=\widetilde{D}_x^{-1}\widetilde{N}_x$ over ${\cal F}({\cal R})$,
$T_x=\bmatrix{N_x^t&D_x^t}^t$ and
$W_x=\bmatrix{\widetilde{N}_x&\widetilde{D}_x}$.  Further for $x=1,2$
and for $I\in{\cal I}$, $J\in{\cal J}$, let $t_{xI}=\det(\Delta_IT_x)$, and $w_{xJ}=\det
(\Delta_JW_x^t)$.  Then the ideals
  $\sum_{I\in{\cal I}}(t_{1I})$, 
  $\sum_{I\in{\cal I}}(t_{2I})$, 
  $\sum_{J\in{\cal J}}(w_{1J})$, and
  $\sum_{J\in{\cal J}}(w_{2J})$ 
are isomorphic to one another as ${\cal R}$-modules.
\end{lemma}
\begin{proof}
We show first~(i) $\sum_{I\in{\cal I}}(t_{1I})\simeq\sum_{I\in{\cal I}}(t_{2I})$ and
then~(ii) $\sum_{I\in{\cal I}}(t_{1I})\simeq\sum_{J\in{\cal J}}(w_{1J})$.  The
isomorphism $\sum_{J\in{\cal J}}(w_{1J})\simeq\sum_{J\in{\cal J}}(w_{2J})$ can be proved
analogously to~(i) and so is omitted.

\noindent (i).  Observe that in the case where
$\bmatrix{N_2^t&D_2^t}^t=\bmatrix{N_1^t&D_1^t}^tX$ holds with some
nonsingular matrix~$X$ over~${\cal R}$, the statement of the lemma obviously
holds.  Hence by considering $\bmatrix{N_2^t&D_2^t}^t\adj(D_2)$ as
$\bmatrix{N_2^t&D_2^t}^t$, we can assume without loss of generality
that~$D_2$ is expressed as $d_2E_m$ with nonzero~$d_2$.  Observe now
that $\bmatrix{N_1^t&D_1^t}^t d_2=\bmatrix{N_2^t&d_2E_m}D_1$ holds.
From this relation and the first observation, we now have~(i).

\noindent
(ii).  It is sufficient to consider the case $P=Nd^{-1}$ with $N\in{\cal R}
^{n\times m}$ and $d\in{\cal R}$ as in~(i).  In the case $P=ND^{-1}$, one can
consider $P=(N\adj(D))\det(D)^{-1}$.

First we define a~bijective mapping~$\tau$ from~${\cal I}$
to~${\cal J}$.  For convenience we decompose~$I$
into~$I_N$ and~$I_d$ as follows
\[
    I_N =\{ i\,|\, i\leq n, i\in I\},{\ \ }
    I_d =\{ i\,|\, i>n,  i\in I\}.
\]
Corresponding to~$I_N$ and~$I_d$,
we define~$J_N$ and~$J_d$ as
\[
    J_N = [1,m]\backslash\{ i-n \,|\, i\in I_d\},{\ \ }
    J_d = \{ i+m\,|\, i\in[1,n]\backslash I_N\}.
\]
We now define the mapping $\tau:{\cal I} \rightarrow{\cal J}$ as
\[
  \tau: I_N\cup I_d \mapsto J_N\cup J_d.
\]
Since~$I_N$ and~$I_d$ can be expressed by~$J_N$
and~$J_d$ as
    $I_N= [1,n]\backslash\{ j-m\,|\, j\in J_d\}$,
    $I_d= \{ j+n\,|\, j\in[1,m]\backslash J_N\}$,
the inverse mapping $\tau^{-1}:{\cal J}\rightarrow{\cal I}$ can be defined naturally.  Hence,
the map~$\tau$ is bijective.

Now let $T=\bmatrix{ N^t & dE_m}^t$ and
$W=\bmatrix{ N & dE_n}$.  By the straightforward
calculation with noting that $dE_m$ and $dE_n$ are
diagonal, we obtain the following relations:
\[%
  \det(\Delta_IT)=\pm\det(\Delta_{\tau(I)}W^t)d^{m-n}.
\]%
Thus $t_{1I}=\pm w_{1 \tau(I)}d^{m-n}$ for all $I\in{\cal I}$.  It follows that
the ideals $\sum_{I\in{\cal I}}(t_{1I})$ and $\sum_{J\in{\cal J}}(w_{1J})$ are
isomorphic to each other.
\qquad
\end{proof}

\begin{note}
\upshape The reduced minors are derived from $t_I$'s and~$w_J$'s in
Definition\ref{D:ReducedMinors}.  Thus $t_I$'s and~$w_J$'s can be
considered more primitive than the reduced minors.  Nevertheless since
we will present in Theorem\,\ref{Th:1stResult} that $t_I$'s and
$w_J$'s (or the ideals~$\mft$ and~$\mfw$ generated by them) have the
capability to state feedback stabilizability over commutative rings,
we here consider that the full-size minor ideal $\mft$ (or the ideal
$\mfw$) is a generalization of the reduced minors.
\end{note}

\section{Feedback Stabilizability in terms of Full-Size Minor Ideal}
\label{S:5}
In this section, we present the necessary and sufficient condition of
the feedback stabilizability over commutative rings in terms of the
full-size minor ideal.

Let us consider the case where the set~${\cal A}$ of the stable causal
transfer functions is not a~unique factorization domain.  Then it is
not sufficient to use the family of reduced minors in order to state
the feedback stabilizability.  To see this, let us consider the result
given by Anantharam in~\cite{bib:anantharam85a}%
\footnote{The 
  author wishes to thank to Dr.~A.~Quadrat (Centre
d'Enseignement et de Recherche en Math\'{e}matiques, Informatique et
Calcul Scientifique, ENPC, France) who introduced him to the paper of
Anantharam\cite{bib:anantharam85a}.}%
.
\begin{example}\label{Ex:Anantharam}
\upshape
In~\cite{bib:anantharam85a}, Anantharam considered the case where
$\Zset[\sqrt{5}\mathrm{i}]$ ($\simeq\Zset[x]/(x^2+5)$) is the set of
stable causal transfer functions, where $\Zset$ is the ring of
integers and~$\mathrm{i}$ the imaginary unit; that is, ${\cal A}
=\Zset[\sqrt{5}\mathrm{i}]$.  The set of all possible transfer
functions is given as the field of fractions of~${\cal A}$; that is, ${\cal F}
=\Qset(\sqrt{5}\mathrm{i})$.  In this case we have multiple
factorizations $2\cdot 3=(1+\sqrt{5}\mathrm{i})(1-\sqrt{5}\mathrm{i})$
over~${\cal A}$, so that~${\cal A}$ is not a~unique factorization domain.
Anantharam in~\cite{bib:anantharam85a} considered the single-input
single-output case and showed that the plant
$p=(1+\sqrt{5}\mathrm{i})/2$ does not have its coprime factorization
over~${\cal A}$ but is stabilizable.

Now let $T=\bmatrix{1+\sqrt{5}\mathrm{i} & 2}^t$.  Since the plant~$p$
is of the single-input single-output ($m=n=1$), we have ${\cal I}=\{\{ 1\},
\{ 2\}\}$.  Thus let $I_1=\{ 1\}$ and $I_2=\{ 2\}$ so that ${\cal I}=\{ I_1,I_2\}
$.  The full-size minors of the matrix~$T$ are
  $t_{I_1}=\det(\Delta_{I_1}T)=1+\sqrt{5}\mathrm{i}$ 
and 
  $t_{I_2}=\det(\Delta_{I_2}T)=2$.  
If Theorem\,\ref{Th:3.4} (or equivalently
Corollary\,\ref{R:19.Jul.99.204321}) could be applied even over
a~general commutative ring, the ideal $(t_{I_1},t_{I_2})$ should be
principal.  However, the ideal $(t_{I_1},t_{I_2})$ is not principal
since~$p$ does not have its coprime factorization.
\QEDBBOX
\end{example}

In order to involve even such an example as a~system over commutative
ring, we extend Theorem\,\ref{Th:3.4}.  Since we cannot
use the reduced minors to state the feedback stabilizability in
general, we alternatively employ the full-size minor ideal~$\mft$
rather than the reduced minors.  The extension is the first main
result of this paper and stated as follows.

\begin{theorem}\label{Th:1stResult}
Let~$P$ be a~causal plant of ${\cal P}^{n\times m}$.  Then the plant~$P$ is
stabilizable if and only if the full-size minor ideal $\mft$ of the
plant~$P$ is projective.  Further when $\mft$ is projective, it is of
rank~$1$.

\end{theorem}

By virtue of Lemma\,\ref{L:4.2}, the above theorem can be also stated
with the ideal~$\mfw$ instead of the full-size minor ideal~$\mft$.

In the case where~${\cal A}$ is a~unique factorization domain, as in
Theorem\,\ref{Th:3.4}, the condition of feedback stabilizability is
that the full-size minor ideal is free.  On the other hand, in
Theorem\,\ref{Th:1stResult}, the condition is that the ideal is
projective.  They are equivalent to each other in the case where~${\cal A}$
is a~unique factorization domain as follows.

\begin{proposition}\label{P:5.3}
Let~${\cal R}$ be a~unique factorization domain.  Then the ideal generated
by finite elements of~${\cal R}$ is projective if and only if it is free.
\end{proposition}

\noindent
This proof will be given after finishing the proof of
Theorem\,\ref{Th:1stResult}.

Now that we have presented the statement of
Theorem\,\ref{Th:1stResult}, the main objective of the remainder of
this section is to carry out the proof of Theorem\,\ref{Th:1stResult}.
To do so, we prepare two main intermediate results.  The first one is
about the existence of right-/left-coprime factorizations of
stabilizable plants over local rings, which will be presented in
Subsection\,\ref{SS:15.Jun.99.103613}.  The other is about the
local-global principle of the feedback stabilizability, which will be
presented in Subsection\,\ref{SS:15.Jun.99.103635}.  Then we will
prove Theorem\,\ref{Th:1stResult}.  After the proof of
Theorem\,\ref{Th:1stResult} we will prove Proposition\,\ref{P:5.3}.
Before finishing this section, we will present the relationship among
the full-size minor ideals of~$P$, $C$, and $H(P,C)$.

\subsection{Right-/Left-Coprime Factorizations over Local Rings}
\label{SS:15.Jun.99.103613}
The following is the first intermediate result of
Theorem\,\ref{Th:1stResult} about the existence of right-/left-coprime
factorizations of stabilizable plants over local rings.

\begin{proposition}\label{P:5.4}
Let~$P$ be a~plant in ${\cal F}^{n\times m}$.  Suppose that~${\cal R}$ is ${\cal A}_\p$ with
a~prime ideal~$\p$ in $\Spec({\cal A})$.  Then the following statements are
equivalent:
\begin{itemize}
\item[\textup{(i)}]  The plant~$P$ is ${\cal R}$-stabilizable.
\item[\textup{(ii)}] There exists a~right-coprime factorization over~${\cal R}$ of~$P$.
\item[\textup{(iii)}] There exists a~left-coprime factorization over~${\cal R}$ of~$P$.
\end{itemize}
\end{proposition}

The proof of this proposition will be presented after giving several
its intermediate results.

We here recall the notion of Hermite used
in~\cite{bib:vidyasagar85a}\footnote{%
It should be noted that this definition of ``Hermite'' is different
from~\cite{bib:kaplansky49a,bib:brmacdonald84a}.
}, which can characterize the existence of both right-/left-coprime
 factorizations of transfer matrices.
\begin{definition}\label{D:5.5}
\textup{\rm (\cite[p.345]{bib:vidyasagar85a})~} Let~${\cal R}$ be a~commutative ring and~$A$ a~matrix over~${\cal R}$ of size $x\times y$ with $x<y$.
Then we say that the matrix~$A$ can be \emph{complemented} if there
exists a~unimodular matrix in $({\cal R})_y$ containing the matrix~$A$ as
a~submatrix.  A~row $\bmatrix{a_1&\cdots&a_y}\in{\cal R}^{1\times y}$ is said to be a
\emph{unimodular row} if $a_1,\ldots, a_y$ together generate~${\cal R}$.  A~commutative ring~${\cal R}$ is said to be \emph{Hermite} if every unimodular
row can be complemented.
\end{definition}

The following result was given in~\cite{bib:vidyasagar85a} provided
that~${\cal R}$ is an integral domain.

\begin{theorem}\label{Th:5.6}
\textup{\rm (cf.~Theorem\,8.1.66 of~\cite{bib:vidyasagar85a})~} Let~$
{\cal R}$ be a~commutative ring.  The following three statements are
equivalent:
\begin{itemize}
\item[\textup{(i)}]~The commutative ring~${\cal R}$ is Hermite.
\item[\textup{(ii)}] 
      If a~matrix over ${\cal F}({\cal R})$ has a~right-coprime factorization
      over~${\cal R}$,
      it has also a~left-coprime factorization over~${\cal R}$.
\item[\textup{(iii)}] 
      If a~matrix over ${\cal F}({\cal R})$ has a~left-coprime factorization
      over~${\cal R}$,
      it has also a~right-coprime factorization over~${\cal R}$.
\end{itemize}
\end{theorem}
The ``integral domain'' version of this theorem was given as
Theorem\,8.1.66 of~\cite{bib:vidyasagar85a}.  Even in the case of
commutative rings, the proof is similar with that of Theorem\,8.1.66
of~\cite{bib:vidyasagar85a} and so is omitted.

The following result is the intermediate result of
Proposition\,\ref{P:5.4}, which makes the result above applicable to
the proof of the proposition.

\begin{lemma}\label{L:5.7}
Any local ring is Hermite. 
\end{lemma}
\begin{proof}
Suppose that~${\cal R}$ is a~local ring and
  $\bmatrix{a_1,\ldots, a_y}\in{\cal R}^{1\times y}$ 
is a~unimodular row.  Thus there exist $b_1,\ldots, b_y\in{\cal R}$ such that
\begin{equation}\label{E:21.May.99.222220}
   a_1b_1+\cdots+a_yb_y=1.
\end{equation}
Since~${\cal R}$ is local, the set of all nonunits is an ideal.  From
(\ref{E:21.May.99.222220}), there exists an~$i$ with $1\leq i\leq y$ such
that~$a_i$ is a~unit.  We assume without loss of generality that~$a_1$
is a~unit.  If $y=1$, then~$a_1$ is a~unit, which can be considered as
a~unimodular matrix of $({\cal R})_1$.  In the following we consider the
case $y>1$.  Then we can construct a~unimodular matrix $U=(u_{ij})\in
({\cal R})_y$:
\[
  u_{ij}=
    \left\{
    \mbox{\begin{tabular}{ll}
            $a_j$       &if $i=1$,\\
            $a_1^{-1}$  &if $i=j=2$,\\
            $1 $        &if $i=j>2$,\\
            $0 $        &otherwise.
          \end{tabular}
         }
    \right.
\]
This~$U$ contains the row $\bmatrix{a_1,\ldots, a_y}$ as a~submatrix and
hence every unimodular row can be complemented.  Therefore~${\cal R}$ is
Hermite.
\qquad
\end{proof}

We prepare one more result which will help us
present a~nonsingular denominator matrix of a~stabilizing controller

\begin{lemma}\label{L:5.9}
Let~${\cal R}$ be a~commutative ring and $\p$ a~prime ideal of~${\cal R}$.
Suppose that there exist matrices~$A$, $B$, $C_1$, $C_2$ over~${\cal R}$
such that the determinant of the following square matrix is in ${\cal R}\backslash
\p$:
\begin{equation}\label{E:L:5.3:1} 
   \bmatrix{ A&C_1 \cr
             B&C_2},
\end{equation}
where the matrix~$A$ is square and the matrices~$A$ and~$B$ have same
number of columns.  Then there exists a~matrix~$R$ over~${\cal R}$ such that
the determinant of the matrix $A+RB$ is in~${\cal R}\backslash\p$.
\end{lemma}

Before starting the proof, it is worth reviewing some easy facts about
a~prime ideal.
\begin{remark}\label{R:5.10}
Suppose that $\p$ is a~prime ideal of~${\cal R}$. (i) If~$a$ is in~${\cal R}\backslash\p$
and expressed as $a=b+c$ with $b,c\in{\cal R}$, then at least one of~$b$
and~$c$ is in~${\cal R}\backslash\p$.
(ii) If~$a$ is in ${\cal R}\backslash\p$ and~$b$ in~$\p$, then the sum $a+b$ is
in ${\cal R}\backslash\p$.
(iii) Every factor in~${\cal R}$ of an element of ${\cal R}\backslash\p$ belongs to ${\cal R}\backslash
\p$ (that is, if $a,b\in{\cal R}$ and $ab\in{\cal R}\backslash\p$, then $a,b\in{\cal R}\backslash\p$).
\end{remark}

\begin{proofof}{Lemma\,\ref{L:5.9}}
This proof mainly follows that of~Lemma\,4.4.21
of~\cite{bib:vidyasagar85a}.

If $\det(A)$ is in ${\cal R}\backslash\p$, then we can select the zero matrix
as~$R$.  Thus we assume in the following that $\det(A)$ is in $\p$.

Since the determinant of~(\ref{E:L:5.3:1}) is in ${\cal R}\backslash\p$, there
exists a~full-size minor of $\bmatrix{A^t & B^t}^t$ in ${\cal R}\backslash\p$ by
Laplace's expansion of~(\ref{E:L:5.3:1}) and by
Remark\,\ref{R:5.10}(i,iii).  Let~$a$ be such a~full-size minor of
  $\bmatrix{A^t & B^t}^t$
having as few rows from~$B$ as possible.  

We here construct a~matrix~$R$ such that
  $\det(A+RB)=\pm a+z$
with~a~$z\in\p$.  Since~$\det(A)\in\nolinebreak\p$, the full-size
minor~$a$ must contain at least one row of~$B$ from the matrix
$\bmatrix{A^t & B^t}^t$.  Suppose that~$a$ is obtained by excluding
the rows $i_1,\ldots, i_k$ of~$A$ and including the rows $j_1,\ldots, j_k$
of~$B$.  Now define $R=(r_{ij})$ by
  $r_{i_1 j_1}=\cdots=r_{i_k j_k}=1$
and $r_{ij}=0$ for all other~$i$,~$j$.  Observe that
$\det(A+RB)$ is expanded in terms of full-size
minors of the matrices $\bmatrix{E & R}$ and
  $\bmatrix{A^t & B^t}^t$
from the factorization
  $A+RB=\bmatrix{E & R}
                       \bmatrix{A^t & B^t}^t$ 
by the Binet-Cauchy formula.  Every minor of $\bmatrix{E & R}$
containing more than~$k$ columns of~$R$ is zero.  By the method of
choosing the rows from 
  $\bmatrix{A^t & B^t}^t$
for the full-size minor~$a$, every full-size minor of
  $\bmatrix{A^t & B^t}^t$
having less than~$k$ rows of~$B$ is in~$\p$.  There is only one
nonzero minor of $\bmatrix{E & R}$ containing exactly~$k$ columns
of~$R$, which is obtained by excluding the columns $i_1,\ldots, i_k$ of
the identity matrix~$E$ and including the columns $j_1,\ldots, j_k$
of~$R$; it is equal to~$\pm 1$.  From the Binet-Cauchy formula the
corresponding minor of $\bmatrix{A^t & B^t}^t$ is~$a$.  As a~result,
$\det(A+RB)$ is given as a {sum} of~$\pm a$ and elements in~$\p$.  By
Remark\,\ref{R:5.10}(ii), the sum is in ${\cal R}\backslash\p$ and so is
$\det(A+RB)$.
\qquad
\end{proofof}

Now that we have the result above, we can prove
Proposition\,\ref{P:5.4}.
\begin{proofof}{Proposition\,\ref{P:5.4}}
Since~${\cal R}$ is local, (ii) and (iii) are equivalent by
Theorem\,\ref{Th:5.6} and Lemma\,\ref{L:5.7}.  Thus we only prove~(i)$
\rightarrow$(ii) and \emph{vice versa}.

\noindent
(i)$\rightarrow$(ii).~~Suppose that~$P$ is ${\cal R}$-stabilizable.  Then the ${\cal R}
$-module ${\cal T}_{P,{\cal R}}$ is projective by Proposition\,2.1
of~\cite{bib:mori97a}.  Further it is free by Corollary\,3.5
of~\cite[Ch.IV]{bib:kunz85a}.  Let~$N$ and~$D$ be matrices over~${\cal R}$
with $P=ND^{-1}\in{\cal F}({\cal R})$.  Then the ${\cal R}$-module $M_r(\bmatrix{N^t &
D^t}^t)$ is free of rank~$m$ since~$D$ is nonsingular over~${\cal R}$.  Let
$v_1, v_2,\ldots, v_m\in{\cal R}^m$ be a~basis of the module
$M_r(\bmatrix{N^t&D^t}^t)$ and~$V$ the matrix of $({\cal R})_m$ whose rows
are $v_1,v_2,\ldots, v_m$.  Then, the matrix $\bmatrix{N^t & D^t}^t$ can
be written in the form
  $\bmatrix{ N^t & D^t}^t=\bmatrix{ N_0^t & D_0^t}^t V$ 
by uniquely choosing the matrices~$N_0$ in ${\cal R}^{n\times m}$ and~$D_0$ in $
({\cal R})_m$.  Because of $\det(D)=\det(D_0V)$, $\det(D_0)$ is a~nonzerodivisor.  It follows that $P=N_0D_0^{-1}$ over ${\cal F}({\cal R})$.  In
the following we show that the matrices~$N_0$ and~$D_0$ are
right-coprime over~${\cal R}$.  Since~$v_1,\ldots, v_m$ belong to
$M_r(\bmatrix{N^t&D^t}^t)$, there exist matrices $\widetilde{Y}$ in $
{\cal R}^{m\times n}$ and $\widetilde{X}$ in $({\cal R})_ m$ such that
    $V=\bmatrix{\widetilde{Y} & \widetilde{X}} \bmatrix{N^t & D^t}^t$.
So we have $V=(\widetilde{Y}N_0+\widetilde{X}D_0)V$.  Since~$V$ is
nonsingular, we obtain
  $\widetilde{Y}N_0+\widetilde{X}D_0=E_m$
over~${\cal R}$.  Thus $(N_0,D_0)$ is a~right-coprime factorization over~${\cal R}
$ of~$P$.

\noindent
(ii)$\rightarrow$(i).~~Suppose that there exists a~right-coprime factorization
over~${\cal R}$ of the plant~$P$; that is, there exist the matrices~$N$,
$D$, $\widetilde{Y}$, $\widetilde{X}$ over~${\cal R}$ with
  $\widetilde{Y}N+\widetilde{X}D=E_m$
and $P=ND^{-1}$.  If $\det(\widetilde{X})$ is a~nonzerodivisor of~${\cal R}
$, it is obvious that $\widetilde{X}^{-1}\widetilde{Y}$ is an ${\cal R}
$-stabilizing controller.  Thus in the following we suppose that
$\det(\widetilde{X})$ is a~zerodivisor of~${\cal R}$.

By the equivalence between~(ii) and~(iii), there also exists
a~left-coprime factorization over~${\cal R}$ of~$P$; that is, there exist
the matrices $\widetilde{N}$, $\widetilde{D}$, $Y$, $X$ over~${\cal R}$ with
  $\widetilde{N}Y+\widetilde{D}X=E_n$
and $P=\widetilde{D}^{-1}\widetilde{N}$.  
Thus we have the following matrix equation:
\begin{equation}\label{E:10.Feb.00.125019}
~~~~~
    \bmatrix{ \widetilde{X} & \widetilde{Y} \cr
              \widetilde{N} &-\widetilde{D}}
    \bmatrix{ D & Y \cr
              N & -X}
    =\bmatrix{ E_m & \widetilde{X}Y-\widetilde{Y}X \cr
               O   & E_n}.
\end{equation}
Observe that the determinant of the right-hand side of the matrix
equation above is in ${\cal R}\backslash{\cal Z}_\p$, where ${\cal Z}_\p$ denotes the
localization of the prime ideal~${\cal Z}$ at~$\p$ (Note that ${\cal Z}_\p$ is
also a~prime ideal of ${\cal R}$).  Hence the determinant of the first
matrix in (\ref{E:10.Feb.00.125019}) is in ${\cal R}\backslash{\cal Z}_\p$ again.
Applying Lemma\,\ref{L:5.9} to the first matrix, we have a~matrix~$R$
over~${\cal R}$ such that the determinant of the matrix
$\widetilde{X}+R\widetilde{N}$ is in~${\cal R}\backslash{\cal Z}_\p$.  Now
$(\widetilde{X}+R\widetilde{N})^{-1}(\widetilde{Y}-R\widetilde{D})$ is
an ${\cal R}$-stabilizing controller.  \qquad
\end{proofof}

\subsection{Local-Global Principle in Stabilizability}
\label{SS:15.Jun.99.103635}
Next we present the local-global principle below about the feedback
stabilizability as the second intermediate result of this section.

\begin{proposition}\label{P:5.8}
Suppose that the plant~$P$ is causal.  Then the following statements
are equivalent:
\begin{romannum}
\item $P$ is stabilizable.
\item $P$ is ${\cal A}_\p$-stabilizable for each prime ideal~$\p$ in
$\Spec({\cal A})$.
\item $P$ is ${\cal A}_\m$-stabilizable for each maximal ideal~$\m$
in $\Max({\cal A})$.
\item For every prime ideal~$\p$ in $\Spec({\cal A})$,
$P$ has either its right- or left-coprime factorization over ${\cal A}_\p$.
\item For every maximal ideal~$\m$ in $\Max({\cal A})$, $P$ has either
its right- or left-coprime factorization over ${\cal A}_\m$.
\end{romannum}
Further, if~$P$ is stabilizable, then there exists a~causal
stabilizing controller of~$P$.
\end{proposition}

Note here that by virtue of Proposition\,\ref{P:5.4}, if~(iv) holds
$\Bigl($if~(v) holds$\Bigr)$, then the plant~$P$ has both
right-/left-coprime factorizations over ${\cal A}_\p$ $\Bigl($over~${\cal A}_
\m\Bigr)$.

We consider that this is a~generalization of Proposition\,2
of~\cite{bib:sule94a} in which the strict causality of the plant is
assumed (see~\cite{bib:sule98a} for the definition of the strict
causality).  On the other hand, we assume only that the plant is
causal. 

Now we begin to prove Proposition\,\ref{P:5.8}.
\begin{proofof}{Proposition\,\ref{P:5.8}}
Since the following implications are obvious:
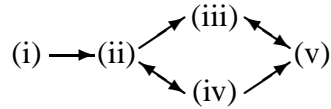
\begin{center}
\unitlength 1mm
\linethickness{0.4pt}
\thicklines
\begin{picture}(48.00,16.00)(0,-8)
\put(  5, 0){\makebox(0,0){(i)}}
\put(  8, 0){\vector(1,0){6.00}}
\put( 17, 0){\makebox(0,0){(ii)}}
\put( 20, 1){\vector(3,2){6.00}}
\put( 20,-1){\vector(3,-2){6.00}}
\put( 20,-1){\vector(-3,2){0.00}}
\put( 30, 5){\makebox(0,0){(iii)}}
\put( 30,-5){\makebox(0,0){(iv)}}
\put( 34, 5){\vector(3,-2){6.00}}
\put( 34, 5){\vector(-3, 2){0.00}}
\put( 34,-5){\vector(3, 2){6.00}}
\put( 43, 0){\makebox(0,0){(v)}}
\end{picture}
\end{center}
by virtue of Proposition\,\ref{P:5.4}, we only show
that~(v) implies~(i).

Suppose that~(v) holds.  Let~$N$, $D$, $\widetilde{N}$, and
$\widetilde{D}$ be matrices over~${\cal A}$ with
$P=ND^{-1}=\widetilde{D}^{-1}\widetilde{N}$ such that~$D$ and
$\widetilde{D}$ are ${\cal Z}$-nonsingular (recall that~$P$ is causal).  By
Proposition\,\ref{P:5.4}, $P$ has both right-/left-coprime
factorizations over ${\cal A}_\m$ with $\m\in\Max({\cal A})$.  As in the proof of
Proposition\,\ref{P:5.4}, for each~$\m$ in $\Max({\cal A})$, there exist
matrices
 $Y_\m$, $X_\m$, $\widetilde{Y}_\m$, $\widetilde{X}_\m$, 
 $N_\m$, $D_\m$, $\widetilde{N}_\m$, $\widetilde{D}_\m$,
 $V_\m$, and $W_\m$
over ${\cal A}_\m$ such that
\begin{eqnarray}
&~~~~~&
   \bmatrix{ N \cr D }
  \!=\!\bmatrix{ N_\m \cr D_\m }V_\m,{\ \ }
   \bmatrix{ \widetilde{N} & \widetilde{D} }
  \!=\!W_\m\bmatrix{ \widetilde{N}_\m & \widetilde{D}_\m },\label{E:28.Dec.99.190551}\\
&~~~~~&
    \widetilde{Y}_\m N_\m+\widetilde{X}_\m D_\m\!=\!E_m,{\ \ }
    \widetilde{N}_\m Y_\m+\widetilde{D}_\m X_\m\!=\!E_n\label{E:28.Dec.99.190810}
\end{eqnarray}
hold over ${\cal A}_\m$.  For each $\m\in\Max({\cal A})$ let $q_\m$ be an
arbitrary but fixed element of~${\cal A}\backslash\m$ such that the six matrices
  $q_\m N_\m\widetilde{Y}_\m$,
  $q_\m N_\m\widetilde{X}_\m$,
  $q_\m D_\m\widetilde{Y}_\m$,
  $q_\m D_\m\widetilde{X}_\m$,
  $q_\m \widetilde{D}_\m$, and
  $q_\m \widetilde{N}_\m$
are over~${\cal A}$.

For a~subset~${\cal B}$ of~${\cal A}$, denote by $\Gamma({\cal B})$ the set of all maximal
ideals~$\m$ of~${\cal A}$ with ${\cal B}\not\subset\m$, that is,
  $\Gamma({\cal B})=\{\m\in\Max({\cal A})\,|\,{\cal B}\not\subset\m\}$.  
Since $q_\m\in{\cal A}\backslash\m$, we have $\m\in \Gamma({\cal A}q_{\m})$.  Thus $\Max({\cal A})= \bigcup_
{\m\in\Max({\cal A})} \Gamma({\cal A}q_{\m})$.  Recall that $\Max({\cal A})$ is compact (see
Theorem\,IV.1 of~\cite{bib:brmacdonald84a}).  Hence there are a~finite
number of $\m_1,\ldots,\m_t$ of maximal ideals such that
  $\Max({\cal A})=\bigcup_{i=1}^t \Gamma({\cal A}q_{\m_i})$. 
It follows that $\Max({\cal A})=\Gamma(\sum_{i=1}^t {\cal A}q_{\m_i})$ and,
consequently, ${\cal A}=\sum_{i=1}^t {\cal A}q_{\m_i}$.  Therefore there exist
$r_1,\ldots, r_t$ in~${\cal A}$ with
  $1=r_1q_{\m_1}+\cdots+r_tq_{\m_t}$.

Next we want to consider that at least one of $\m_1,\ldots,\m_t$
contains~${\cal Z}$.  In the case where every $\m_i$ in $\m_1,\ldots,\m_t$ does
not contain~${\cal Z}$, we reconstruct~$t$, $r_i$'s, and $q_{\m_i}$'s as
follows.  We first pick an $\m_{t+1}\in\Max({\cal A})$ with $\m_{t+1}\supset{\cal Z}$.
Then we let~$r_i$ be $(1-q_{\m_{t+1}})r_i$ for $1\leq i\leq t$ and
$r_{t+1}=1$.  We now let $t:=t+1$.  Then we have again
  $1=r_1q_{\m_1}+\cdots+r_tq_{\m_t}$
and, in this case, $\m_t\supset{\cal Z}$.  Hence we can assume without loss of
generality that at least one of $\m_1,\ldots,\m_t$, say $\m_1$,
contains~${\cal Z}$.

Observe then that the following equality holds:
\begin{equation}\label{E:23.Nov.99.144001}
\mbox{$\begin{array}{r}
   1=(r_1q_{\m_1}+r_1-1)q_{\m_1}+(r_2q_{\m_1}+r_2)q_{\m_2}~~~~\\
                             +\cdots+(r_tq_{\m_1}+r_t)q_{\m_t}.
\end{array}$}
\end{equation}
At least one of $r_1q_{\m_1}+r_1-1$ and~$r_1$ must be in ${\cal A}\backslash{\cal Z}
$.  Thus in the case $r_1\in{\cal Z}$, we can reassign~$r_i$'s as in
(\ref{E:23.Nov.99.144001}), so that~$r_1$ is in ${\cal A}\backslash{\cal Z}$.  Therefore
we can assume without loss of generality that $r_1q_{\m_1}\in{\cal A}\backslash{\cal Z}$.

Consider here the following matrix
\begin{equation}\label{E:09.Jun.99.123818}
\makebox[0cm]{$
\hspace*{10mm}
   \bmatrix{\!   E_n-\sum_{i=1}^t r_iq_{\m_i} N_{\m_i}\widetilde{Y}_{\m_i} &
                -\sum_{i=1}^t r_iq_{\m_i} N_{\m_i}\widetilde{X}_{\m_i} \cr
            \!   \sum_{i=1}^t r_iq_{\m_i} D_{\m_i}\widetilde{Y}_{\m_i} &
                 \sum_{i=1}^t r_iq_{\m_i} D_{\m_i}\widetilde{X}_{\m_i}},
$}
\end{equation}
which is over~${\cal A}$.  For short we partition (\ref{E:09.Jun.99.123818}) as
\[
  \bmatrix{ H_{11} & H_{12} \cr
            H_{21} & H_{22}}.
\]
In the case where $H_{22}$ is ${\cal Z}$-nonsingular, letting
$C=H_{22}^{-1}H_{21}\in{\cal P}^{m\times n}$ we can check that $H(P,C)$ is equal
to (\ref{E:09.Jun.99.123818}), which implies that~$P$ is stabilized
by~$C$.
Hence in the rest of this proof we show that if $H_{22}$ is ${\cal Z}
$-singular, then $H_{22}$ can be made ${\cal Z}$-nonsingular by reassigning
$\widetilde{X}_{\m_i}$ and $\widetilde{Y}_{\m_i}$ for an~$i$.

First we show the ${\cal Z}$-nonsingularity of the matrices
$r_1q_{\m_1}D_{\m_1}$ and $r_1q_{\m_1}\widetilde{D}_{\m_1}$.  Since
$r_1q_{\m_1}\in{\cal A}\backslash{\cal Z}$, we have $\det(r_1q_{\m_1}D)\in{\cal A}\backslash{\cal Z}$.  From the
first matrix equation of (\ref{E:28.Dec.99.190551}), we have
$\det(r_1q_{\m_1}D)=\det(r_1q_{\m_1}D_{\m_1})\det(V_{\m_1})$.  Hence
$r_1q_{\m_1}D_{\m_1}$ is ${\cal Z}$-nonsingular.  Analogously, from the
second matrix equation of (\ref{E:28.Dec.99.190551}),
$r_1q_{\m_1}\widetilde{D}_{\m_1}$ is ${\cal Z}$-nonsingular.

Next consider the following matrix equation over~${\cal A}$:
\begin{eqnarray}
&
\makebox[0cm]{$
    \bmatrix{
      \sum_{i=1}^t r_iq_{\m_i} D_{\m_i}\widetilde{X}_{\m_i} \!\!&\!\!
      \sum_{i=1}^t r_iq_{\m_i} D_{\m_i}\widetilde{Y}_{\m_i} \cr
      \!\! -r_1q_{\m_1} \det(r_1q_{\m_1}D_{\m_1})\widetilde{N}_{\m_1}  \!\!&\!\!
       r_1q_{\m_1} \det(r_1q_{\m_1}D_{\m_1})\widetilde{D}_{\m_1} }\times
$}&
\nonumber\\
\makebox[0cm]{\hspace*{3em}$
    \bmatrix{ D & O \cr
              N & E_n
     }
    =
    \bmatrix{ D & \sum_{i=1}^t r_iq_{\m_i} D_{\m_i}\widetilde{Y}_{\m_i} \cr
              O & r_1q_{\m_1}\det(r_1q_{\m_1}D_{\m_1})\widetilde{D}_{\m_1} }.
$}&\label{E:08.Feb.97.220842}
\end{eqnarray}
Since the matrices~$D$, $r_1q_{\m_1}D_{\m_1}$, and
$r_1q_{\m_1}\widetilde{D}_{\m_1}$ are ${\cal Z}$-nonsingular, so is the
right-hand side of (\ref{E:08.Feb.97.220842}).  Thus the first matrix
of (\ref{E:08.Feb.97.220842}) is also ${\cal Z}$-nonsingular.  By
Lemma\,\ref{L:5.9} and the first matrix of (\ref{E:08.Feb.97.220842}),
there exists a~matrix $R_{\m_1}'$ of ${\cal A}^{m\times n}$ such that the
following matrix is ${\cal Z}$-nonsingular:
\[%
  \sum_{i=1}^t r_iq_{\m_i} D_{\m_i}\widetilde{X}_{\m_i}
  -r_1q_{\m_1} \det(r_1q_{\m_1}D_{\m_1})
                                  R_{\m_1}'\widetilde{N}_{\m_1}.
\]%
Now let $R_{\m_1}$ be
  $r_1q_{\m_1}\adj(r_1q_{\m_1}D_{\m_1})R_{\m_1}'$.
Further we let $\widetilde{X}_{\m_1}$ be the matrix
  $\widetilde{X}_{\m_1}-R_{\m_1}\widetilde{N}_{\m_1}$
and $\widetilde{Y}_{\m_1}$  the matrix
  $\widetilde{Y}_{\m_1}+R_{\m_1}\widetilde{D}_{\m_1}$,
which are consistent with (\ref{E:28.Dec.99.190810}).  Thus we can now
consider without loss of generality that the matrix $ \sum_ {i=1}^t
r_iq_{\m_i} D_{\m_i}\widetilde{X}_{\m_i}$ is ${\cal Z}$-nonsingular and so
is~$H_{22}$.
\qquad
\end{proofof}

\subsection{Proof of Theorem\,\ref{Th:1stResult}}
\label{SS:14.Jul.99.192111}
Before proving Theorem\,\ref{Th:1stResult}, we should prepare a~small
result.

\begin{lemma}\label{L:5.11}
Let $a\in{\cal A}$ and $\p\in\Spec({\cal A})$.  Then $(a)_\p$ and $(a/1)$ are
isomorphic to each other as ${\cal A}_\p$-modules, where $(a)_\p$ denotes
the localization, at~$\p$, of the principal ideal generated by $a\in{\cal A}$
and $(a/1)$ the principal ideal generated by $a/1\in{\cal A}_\p$.
\end{lemma}

The proof of the lemma is elementary and is omitted.

Now we start to prove the first result of this paper.
\begin{proofof}{Theorem\,\ref{Th:1stResult}}
We show first the ``Only If'' part and then the ``If'' part.

\noindent
(Only If).  Suppose that~$P$ is stabilizable.  Then by
Proposition\,\ref{P:5.8}, for every prime ideal~$\p$ in $\Spec({\cal A})$,
$P$ is ${\cal A}_\p$-stabilizable.  By Proposition\,\ref{P:5.4}, $P$ has
both its right-/left-coprime factorizations over~${\cal A}_\p$.  Suppose
that $\widetilde{Y}_{\p}N_{\p}+\widetilde{X}_{\p}D_{\p}=E_m$ holds
over ${\cal A}_\p$ with $P=N_{\p}D_{\p}^{-1}$, where the matrices $N_\p$,
$D_\p$, $\widetilde{Y}_\p$, and $\widetilde{X}_\p$ are over ${\cal A}_\p$.
Then let $T_\p=\bmatrix{N_\p^t&D_\p^t}^t$.  By Binet-Cauchy formula we
have $\sum_{I\in{\cal I}}(\det(\Delta_IT_\p))={\cal A}_\p$.  Thus by virtue of
Lemmas\,\ref{L:4.2} and~\ref{L:5.11}, the ideal $\mft_\p$ is free
(recall that $\mft_\p$ denotes the localization of the full-size minor
ideal $\mft$ at $\p$), which is also finitely generated.  This holds
for every prime ideal~$\p$.  From Theorem IV.32
of~\cite{bib:brmacdonald84a}, the full-size minor ideal~$\mft$ is
projective.

\noindent
(If).  Suppose that the full-size minor ideal $\mft$ is projective.
Let~$\p$ be a~prime ideal in $\Spec({\cal A})$.  Then $\mft_\p$ is free by
Theorem IV.32 of~\cite{bib:brmacdonald84a} again.  Thus there
exist~$g$, $a_I$, and~$r_I$ in ${\cal A}_\p$ with $g=\sum_ {I\in{\cal I}} r_It_I$
and $t_I=a_Ig$ for every $I\in{\cal I}$.  Since $g= \sum_{I\in{\cal I}} r_Ia_Ig$
and~$g$ is a~nonzerodivisor, we have $\sum_{I\in{\cal I}} r_Ia_I=1$.  Recall
here that~${\cal A}_\p$ is local.  Hence the set of all nonunits in ${\cal A}_\p$
is an ideal.  Thus there exists $I_0\in{\cal I}$ such that $r_{I_0}a_{I_0}$
is a~unit of ${\cal A}_\p$.  This implies that $a_{I_0}$ is a~unit of ${\cal A}_
\p$ and further that every~$t_I$ has a~factor $t_{I_0}$ over~${\cal A}_\p$
(that is, $t_{I_0}$ and~$g$ are associate).  Now let $T'=T\adj (\Delta_
{I_0}T)$ and $t_I'=\det (\Delta_IT')$ for every $I\in{\cal I}$.  Then
$t_I'=t_I\det(\adj (\Delta_ {I_0}T))$ and $\Delta_ {I_0}T'=t_{I_0}E_m$ hold.
Since $\det(\adj(\Delta_ {I_0}T))=t_{I_0}^{m-1}$, every~$t_I'$ has
a~common factor $t_{I_0}^m$.

Suppose that~$i$ is an integer with $i\not\in I_0$ and $1\leq i\leq m+n$.  Suppose
further that $i_{01},\ldots, i_{0m}$ are elements in~$I_0$ with ascending
order.  Now let $I=\{ i, i_{01},i_{02},\ldots, i_{0\,k-1},$
$i_{0\,k+1},\ldots, i_{0m}\}$.  Then~$t_I$ is expressed as $\pm
t_{ik}t_{I_0}^{m-1}$ where $t_{ik}$ is the $(i,k)$-entry of the
matrix~$T'$.  Since~$t_I'$ has a~factor~$t_{I_0}^m$, $t_{ik}$ has
a~factor~$t_{I_0}$.  This fact holds for all~$i$ between $1\leq i\leq m+n$
but $i\not\in I_0$.  As a~result, $t_{I_0}$ is a~common factor of all
entries of~$T'$.

Let $T''=T'/t_{I_0}$ over ${\cal A}_\p$.  Since $\Delta_{I_0}T''$ is the
identity matrix, the matrix $\Delta_{I_0}$ itself is a~left inverse of
$T''$.  Let $\widetilde{Y}_{I_0}$ and $\widetilde{X}_{I_0}$ be
matrices with $\bmatrix{\widetilde{Y}_{I_0}&\widetilde{X}_{I_0}}=\Delta_
{I_0}$.  Further we let~$N_{I_0}$ and~$D_{I_0}$ be matrices over ${\cal A}_
\p$ with $T''=\bmatrix{N_{I_0}^t&D_{I_0}^t}^t$.  Then we obtain
   $\widetilde{Y}_{I_0}N_{I_0}+\widetilde{X}_{I_0}D_{I_0}=E_m$
over ${\cal A}_\p$, which is a~right-coprime factorization over ${\cal A}_\p$ of
the plant~$P$.  Therefore by Proposition\,\ref{P:5.8}, $P$
is stabilizable.
\qquad
\end{proofof}

\subsection{Proof of Proposition\,\ref{P:5.3}}
Now we prove Proposition\,\ref{P:5.3}.  We first prepare the following
local-global principle on ideals.
\begin{lemma}\label{L:5.12}
Let~${\cal R}$ be a~commutative ring.  Let $\mfa_1,\ldots,\mfa_k$ be ideals
of~${\cal R}$.  Then the following statements are equivalent:
\begin{itemize}
\item[\textup{(i)}]   $\mfa_1+\cdots+\mfa_k={\cal R}$.
\item[\textup{(ii)}]  $\mfa_{1\p}+\cdots+\mfa_{k\p}={\cal R}_\p$
                    for all prime ideal $\p\in\Spec({\cal A})$.
\item[\textup{(iii)}] $\mfa_{1\m}+\cdots+\mfa_{k\m}={\cal R}_\m$
                    for all maximal ideal $\m\in\Max({\cal A})$.
\end{itemize}
\end{lemma}
\begin{proof}
It is obvious that~(i) implies~(ii) and~(ii) implies~(iii).  Hence we
only show that~(iii) implies~(i).

\noindent
(iii)$\rightarrow$(i).  Suppose that~(iii) holds.  Let~$\m$ be a~maximal ideal
of~${\cal A}$.  Since ${\cal R}_\m$ is local, the set of all nonunits in ${\cal R}_\m$
is an ideal.  Hence there exists an $i_{\m}$ with $1\leq i_\m\leq k$ such
that $\mfa_{i_{\m}\m}={\cal R}_\m$.  Thus there exists $s_\m$ in ${\cal R}\backslash\m$
such that $s_\m\in\mfa_{i_{\m}}$.

Recalling the proof of Proposition\,\ref{P:5.8}, we have a~finite
number of $\m_1,\ldots,\m_t$ in $\Max({\cal R})$ and $r_1,\ldots, r_t\in{\cal R}$ such
that $1=r_1s_{\m_1}+\cdots+r_ts_{\m_t}$ over~${\cal R}$.   
For every $l=1$ to~$t$, $r_ls_{\m_l}$ is an element of $\mfa_i$ with
$i=i_{\m_l}$.  Therefore we have~(i).  \qquad
\end{proof}

\begin{proofof}{Proposition\,\ref{P:5.3}}
Suppose that~${\cal R}$ is a~unique factorization domain.  Since the ``If''
part is obvious, we prove only the ``Only If'' part.  

\noindent 
(Only If).  Let $a_1,\ldots, a_k$ be in~${\cal R}$.  Suppose that
$(a_1,\ldots, a_k)$ is projective.  If all $a_1,\ldots, a_k$ are zero, the
proof is obvious.  Thus in the following we suppose that at least one
of $a_1,\ldots, a_k$ is nonzero.  Since~${\cal R}$ is a~unique factorization
domain, there exists a~nonzero greatest common factor of $a_i$'s,
denoted by~$g$.  Thus there exist $b_i$'s in~${\cal A}$ with $b_ig=a_i$.
Then $(b_1,\ldots, b_k)$ is projective again.  For any prime ideal~$\p$ in
$\Spec({\cal R})$, $(b_1,\ldots, b_k)_\p$ is free of rank~$1$.  Since there is
no nonunit common factor among $b_i$'s over~${\cal R}$, $(b_1,\ldots, b_k)_\p={\cal R}_
\p$.  By Lemma\,\ref{L:5.12}, $(b_1,\ldots, b_k)={\cal R}$.  Hence
$(a_1,\ldots, a_k)=(g)$, which is free.
\qquad
\end{proofof}

\subsection{Full-Size Minor Ideals of~$P$, $C$, and $H(P,C)$}
Now that we have obtained Theorem\,\ref{Th:1stResult}, we know that
the projectivity of the full-size minor ideal of the plant connects
with the feedback stabilizability of the plant.  Since~$P$, $C$, and
$H(P,C)$ are transfer matrices over~${\cal F}$, we can define the full-size
minor ideals of~$C$ and $H(P,C)$ analogously to that of~$P$.

We present here the relationship among the
full-size minor ideals of~$P$, $C$, and $H(P,C)$.

\begin{proposition}\label{P:5.13}
Let $\mft_P$, $\mft_C$, $\mft_{H(P,C)}$ be the full-size minor ideals
of~$P$, $C$, and $H(P,C)$, respectively.  Then $\mft_{H(P,C)}$ is
isomorphic (as an ${\cal A}$-module) to the ideal generated by $t_1t_2$'s
for all $t_1\in\mft_P$ and all $t_2\in\mft_C$.
\end{proposition}

This proposition holds even if~$C$ is not a~stabilizing controller
of~$P$.  Before proving this proposition, we present a~preliminary
lemma.
\begin{lemma}\label{L:5.14}
Let~$A$ and~$B$ are matrices over~${\cal R}$ such that $B=UA$, where~$U$ is
a~unimodular matrix over~${\cal R}$.  Then the ideal generated by the
full-size minors of~$A$ is equal to that of~$B$.
\end{lemma}

The proof of this lemma is straightforward and  omitted.

\begin{proofof}{Proposition\,\ref{P:5.13}}
By virtue of Lemma\,\ref{L:4.2}, we suppose without loss of generality
that~$N$ and~$N_c$ are matrices over~${\cal A}$ and~$d$ and~$d_c$ in~${\cal A}$
with $P=Nd^{-1}$ and $C=N_cd_c^{-1}$.  Let~$A$ and~$B$ be the
following matrices:
\begin{eqnarray*}
&&
    A=\bmatrix{ N_c & O \cr
                d_cE_n & O \cr
                O & N \cr
                O & dE_m},{\ \ }
    B=\bmatrix{ Q \cr
                S }, \mbox{ where}\\
&&
    Q=\bmatrix{ d_cE_n & N \cr
                -N_c & dE_m},{\ \ }
    S=\bmatrix{ d_cE_n & O \cr
                 O & dE_m}.
\end{eqnarray*}
Then we can see that there exists a~unimodular matrix~$U$ with $B=UA$
and that $H(P,C)=SQ^{-1}$.  Let $\mfa$ be the ideal generated by the
full-size minors of~$A$ and $\mft_{P,C}$ be the ideal generated by
$t_1t_2$'s for all $t_1\in\mft_P$ and all $t_2\in\mft_C$.  Then by
Lemma\,\ref{L:5.14}, $\mft_{H(P,C)}$ is isomorphic to $\mfa$ as ${\cal A}
$-modules.  Also by Binet-Cauchy formula, $\mfa\simeq\mft_{P,C}$.  Hence
we obtain $\mft_{H(P,C)}\simeq\mft_{P,C}$.
\qquad
\end{proofof}

\section{Stabilizability in terms of Coprimeness of Quotient Ideals}
\label{S:6}

In this section, we present one more necessary and sufficient
condition of the feedback stabilizability which is given in terms of
quotient ideals.

\begin{theorem}\label{Th:6.1}
Let~$P$ be a~causal plant of ${\cal P}^{n\times m}$.  Then the plant~$P$ is
stabilizable if and only if the ideal
\begin{equation}\label{E:6.1}
    \sum_{I\in{\cal I}} ((t_I):\mft)
\end{equation}
is equal to~${\cal A}$.
\end{theorem}

The ideal of (\ref{E:6.1}) will be considered as another
generalization of the reduced minors.  This will be presented later as
Proposition\,\ref{P:6.7}.

We note that the result above can be considered as a~generalization of
Theorem\,2.1.1 in~\cite{bib:shankar92a} given by Shankar and Sule as
well as a~generalization of Theorem\,\ref{Th:3.4}.  They considered
the single-input single-output case.  In Theorem\,2.1.1
of~\cite{bib:shankar92a}, they stated the feedback stabilizability of
the given plant in terms of the coprimeness of the ideal quotients as
(\ref{E:6.1}).  As a~result, Theorem\,\ref{Th:6.1} can be considered
as a~multi-input multi-output version of Theorem\,2.1.1
of~\cite{bib:shankar92a}.

In order to prove Theorem\,\ref{Th:6.1}, we prepare a~relationship
between projective modules and quotient ideals as follows.

\begin{theorem}\label{Th:6.2}
Let~${\cal R}$ be a~commutative ring and $a_1,\ldots, a_k\in{\cal R}$.  Then
$(a_1,\ldots, a_k)$, that is, the ideal generated by $a_1,\ldots, a_k$ is
projective if and only if the following equation holds:
\begin{equation}\label{E(Th:6.2):1}
    \sum_{i=1}^k ((a_i):(a_1,\ldots, a_k))={\cal R}.
\end{equation}
\end{theorem}

Once we obtain Theorem\,\ref{Th:6.2}, the proof of
Theorem\,\ref{Th:6.1} is directly obtained from
Theorems\,\ref{Th:1stResult} and~\ref{Th:6.2}.  Thus we will present
only the proof of Theorem\,\ref{Th:6.2}, which will be given after
showing intermediate results (Lemmas\,\ref{L:6.3} and~\ref{L:6.4}).

\begin{lemma}\label{L:6.3}
Let~${\cal R}$ be a~commutative ring and $a_1,\ldots, a_k\in{\cal R}$.  If
$(a_1,\ldots, a_k)$ is free, then (\ref{E(Th:6.2):1}) holds.
\end{lemma}
\begin{proof}
As in the proof of Proposition\,\ref{P:5.3}, if all $a_1,\ldots, a_k$ are
zero, the proof is obvious.  Thus in the following we assume that at
least one of $a_1,\ldots, a_k$ is nonzero.  Then there exist a~nonzero $g$
in~${\cal R}$ and~$b_i$ in~${\cal R}$ for $i=1$ to~$k$ such that
$(g)=(a_1,\ldots, a_k)$ and $a_i=b_ig$.  Thus there exist $r_i\in{\cal R}$ for
$i=1$ to~$k$ with $g=r_1a_1+\cdots+r_ka_k$.  If~$g$ was a~zerodivisor, the
principal ideal $(g)$ could not be free.  Hence~$g$ is a~nonzerodivisor.  Now we have
\begin{equation}\label{E:03.Jun.99.131228}
  r_1b_1+\cdots+r_kb_k=1.
\end{equation}
Since $b_i(a_1,\ldots, a_k)\subset(a_i)$ for all~$i$, we have $b_i\in
((a_i):(a_1,\ldots, a_k))$.  It follows from~(\ref{E:03.Jun.99.131228})
that we now have (\ref{E(Th:6.2):1}).  
\qquad
\end{proof}

\begin{lemma}\label{L:6.4}
Let~${\cal R}$ be a~commutative ring, $\mfa,\mfb$ ideals of~${\cal R}$, and~$\p$
a~prime ideal of~${\cal R}$.  Denote by $(\mfa:\mfb)_\p$ the localization of
the quotient ideal $(\mfa:\mfb)$ at~$\p$.  Further let
$(\mfa_\p:\mfb_\p)$ be the quotient ideal of ${\cal R}_\p$, where $\mfa_\p$
and $\mfb_\p$ are localizations of ideals~$\mfa$ and~$\mfb$ at~$\p$,
respectively.  Then $(\mfa:\mfb)_\p=(\mfa_\p:\mfb_\p)$ holds.
\end{lemma}

Now we are in a~position to prove Theorem\,\ref{Th:6.2}.
\begin{proofof}{Theorem\,\ref{Th:6.2}}
By the same reason as in the proofs of Proposition\,\ref{P:5.3} and
Lemma\,\ref{L:6.3}, we assume that at least one of $a_1,\ldots, a_k$ is
nonzero.

\noindent
(If).  Suppose that (\ref{E(Th:6.2):1}) holds.  Then there exist $x_i
\in((a_i):(a_1,\ldots, a_k))$ for $i=1$ to~$k$ such that $1= \sum_ {i=1}^k
x_i$.  By appropriate changes of $a_1,\ldots, a_k$, we assume without loss
of generality that all $x_1,\ldots, x_{k'}$ are nonzero with $1\leq k'\leq k$
and all $x_{k'+1},\ldots, x_k$ are zero subject to $k'<k$.  Observe that
for each~$i$ between~$1$ and~$k'$, $(a_1,\ldots, a_k)_{x_i}=(a_i)_{x_i}$
over ${\cal A}_{x_i}$, where $(a_1,\ldots, a_k)_{x_i}$ and $(a_i)_{x_i}$ denote
the localizations of $(a_1,\ldots, a_k)$ and $(a_i)$ at~$x_i$,
respectively.  Hence for each~$i$ between~$1$ and~$k'$,
$(a_1,\ldots, a_k)_{x_i}$ is free over ${\cal A}_ {x_i}$.  Therefore by Theorem
IV.32 of~\cite{bib:brmacdonald84a}, $(a_1,\ldots, a_k)$ is projective as $
{\cal R}$-module.

\noindent 
(Only If).  Suppose that $(a_1,\ldots, a_k)$ is projective.  Then again by
Theorem IV.32 of~\cite{bib:brmacdonald84a}, for each~$\p$ in $\Spec({\cal R})
$, $(a_1,\ldots, a_k)_\p$ is free over ${\cal R}_\p$.  By
Lemma\,\ref{L:6.3}, we have
\begin{equation}\label{E:03.Jun.99.162501}
    \sum_{i=1}^k ((a_i)_\p:(a_1,\ldots, a_k)_\p)={\cal R}_\p
\end{equation}
for each~$\p$ in $\Spec({\cal R})$.  Then (\ref{E:03.Jun.99.162501}) can be
rewritten as follows by Lemma\,\ref{L:6.4}:
\begin{equation}\label{E:03.Jun.99.162501:2}
    \sum_{i=1}^k ((a_i):(a_1,\ldots, a_k))_\p={\cal R}_\p.
\end{equation}
Since this holds for every $\p$ in $\Spec({\cal A})$, applying
Lemma\,\ref{L:5.12} to (\ref{E:03.Jun.99.162501:2}) we obtain
(\ref{E(Th:6.2):1}).
\qquad
\end{proofof}

We now connect the reduced minors with the quotient ideal of
(\ref{E:6.1}) provided that~${\cal A}$ is a~unique factorization domain.

\begin{proposition}\label{P:6.7}
Suppose that~${\cal A}$ is a~unique factorization domain.  Let~$a_I$ denote
the reduced minor of the matrix~$T$ with respect to $I\in{\cal I}$.  Then
$(a_I)=((t_I):\mft)$ holds for every $I\in{\cal I}$.
\end{proposition}
\begin{proof}
We first show~(i) $(a_I)\subset((t_I):\mft)$ and then~(ii) the opposite
inclusion.

\noindent
(i).  For every $I'\in{\cal I}$, $a_It_{I'}=a_{I'}t_I$ holds, which implies
that $a_I\in((t_I):(t_{I'}))$.  Hence $a_I\in((t_I):\mft)$.

\noindent
(ii).  Suppose that~$\lambda_I$ is an element of the quotient ideal
$((t_I):\mft)$.  Then for every $I'\in{\cal I}$, there exists $\nu_{I'}\in{\cal A}$
such that $\lambda_ It_{I'}=\nu_{I'}t_I$ holds and so $\lambda_Ia_{I'}=\nu_
{I'}a_I$.  Since this equality holds for every $I'\in{\cal I}$, $\lambda_I$ has
a~factor~$a_I$.  Hence $\lambda_I\in(a_I)$.
\qquad
\end{proof}

From the result above, the reduced minor of the matrix~$T$ with
respect to $I\in{\cal I}$ is equal to the quotient ideal $((t_I):\mft)$ up to
a~unit multiple of~${\cal A}$ provided that ${\cal A}$ is a unique factorization
domain.

Now that we have shown a new criterion (\ref{E:6.1}) of the feedback
stabilizability, in the following we present the relationship between
generalized elementary factors and (\ref{E:6.1}) by using radicals of
ideals.
\begin{theorem}
\label{Th:15.Jun.99.114014}
Let $\Lambda_{P\!I}$ denote the generalized elementary factor of the
plant~$P$ with respect to~$I$ in~${\cal I}$.  Then the radical of $\Lambda_
{P\!I}$ is equal to the radical of $((t_I):\mft)$.
\end{theorem}

Before proving this result, we present an analogous result of
Lemma\,\ref{L:6.4}.
\begin{lemma}\label{L:03.Jun.99.162445:x}
Let~${\cal R}$ be a~commutative ring, $\mfa,\mfb$ ideals of~${\cal R}$, and $f\in{\cal R}
$.  Denote by $(\mfa:\mfb)_f$ the localization of the quotient ideal
$(\mfa:\mfb)$ at~$f$.  Further let $(\mfa_f:\mfb_f)$ be the quotient
ideal of~${\cal R}_ f$, where $\mfa_f$ and $\mfb_f$ are localizations of
principal ideals $\mfa$ and $\mfb$ at~$f$, respectively.  Then
$(\mfa:\mfb)_f=(\mfa_f:\mfb_f)$ holds.
\end{lemma}

Analogously to Lemma\,\ref{L:6.4}, the proof of this lemma is omitted.

\begin{proofof}{Theorem\,\ref{Th:15.Jun.99.114014}}
Let $I$ be fixed.  We first show~(i) $\Lambda_{P\!I}\subset\sqrt{((t_I):\mft)}$
and then~(ii) $\sqrt{\Lambda_{P\!I}}\supset((t_I):\mft)$.  They are sufficient
to prove this theorem.

\noindent
(i).  Let~$\lambda$ be an arbitrary but fixed element of $\Lambda_{P\!I}$.  Then
there exists a~matrix~$K$ over~${\cal A}$ with $\lambda T =K \Delta_I T$.  Then for
every $I'\in{\cal I}$, we have $\lambda\Delta_{I'}T =\Delta_{I'}K \Delta_I T$, so that $\lambda
^mt_{I'}=\det(\Delta_ {I'}K)t_I$.  This implies $\lambda^m\in((t_I):(t_{I'}))$.
Hence we have
  $\lambda^m\in\bigcap_{I'\in{\cal I}}((t_I):(t_{I'}))=
             ((t_I):\sum_{I'\in{\cal I}} (t_{I'}))$.

\noindent
(ii).  Let~$\lambda$ be an arbitrary but fixed element of $((t_I):\mft)$.
Then
   $((t_I):\mft)_{\lambda} ={\cal A}_{\lambda}$
and hence
   $((t_I)_{\lambda}:\mft_{\lambda}) ={\cal A}_{\lambda}$
by Lemma\,\ref{L:03.Jun.99.162445:x}.  This implies that $(t_I)_{\lambda}
=\mft_{\lambda}$ and further that every full-size minor of~$T$ has
a~factor~$t_I$ over~${\cal A}_{\lambda}$.  Since~$t_I$ is a~factor of $\det(D)$, it
is a~nonzerodivisor of~${\cal A}_{\lambda}$.  Now let $T'=T(\adj (\Delta_ IT))$ and
$t'_{I'}=\det (\Delta_ {I'}T')$ for every $I'\in{\cal I}$.  Then
$t'_{I'}=t_{I'}\det(\adj (\Delta_IT))$ and $\Delta_IT'=t_IE_m$ hold.  Since
$\det(\adj(\Delta_IT))=t_I^{m-1}$, every~$t_{I'}'$ has a~common
factor~$t_I^m$.

Analogously to the proof of Theorem\,\ref{Th:1stResult}, we can
show that every entry of~$T'$ has a~factor~$t_I$.  Let $T''=T'/t_I$
over~${\cal A}_{\lambda}$.  Then $T=T''\Delta_IT$ holds over~${\cal A}_{\lambda}$.  Hence there
exists an integer~$\omega$ such that $\lambda_I^{\omega}T''$ can be considered over~$
{\cal A}$ and further $\lambda_I^{\omega}T=\lambda_I^{\omega}T''\Delta_IT$ holds over~${\cal A}$.  Now
letting $K=\lambda_I^{\omega}T''\Delta_I$, we have that~$\lambda_I^{\omega}$ is an element of $
\Lambda_{P\!I}$ and hence $\lambda_I\in\sqrt{\Lambda_ {P\!I}}$.  \qquad
\end{proofof}

In the case where~${\cal A}$ is a~unique factorization domain, we obtain the
following result which connects Theorems\,\ref{Th:3.5}
and~\ref{Th:3.4}.
\begin{theorem}
\label{Th:15.Jun.99.161608}
Suppose that~${\cal A}$ is a~unique factorization domain.  Let~$P$ be
a~causal plant and~$I$ in~${\cal I}$.  Then the radical of the elementary
factor of the matrix~$T$ with respect to~$I$ is equal to the radical
of the reduced minor of~$T$ with respect to~$I$ up to a~unit multiple.
\end{theorem}
\begin{proof}
Let~$f_I$ denote the elementary factor of the matrix~$T$ with respect
to~$I$.  Also let~$a_I$ denote the reduced minor of~$T$ with respect
to~$I$.

In the case where~${\cal A}$ is a~unique factorization domain, the
generalized elementary factor of the plant~$P$ with respect to~$I$ is
equal to the principal ideal $(f_I)$.  Thus, by
Theorem\,\ref{Th:15.Jun.99.114014},
$\sqrt{(f_I)}=\sqrt{((t_I):\mft)}$.  By virtue of
Proposition\,\ref{P:6.7}, we have $\sqrt{(f_I)}=\sqrt{(a_I)}$.
\qquad%
\end{proof}

\section{Concluding Remarks}
We have presented two generalization of the reduced minors.  One is
the full-size minor ideal.  Its projectivity is a criterion of the
feedback stabilizability(Theorem\,\ref{Th:1stResult}).  The other is
quotient ideals in (\ref{E(Th:6.2):1}).  Their coprimeness is a
criterion of the feedback stabilizability(Theorem\,\ref{Th:6.1}).

\end{document}